\documentclass[reqno]{amsart}
\usepackage{amssymb,amsmath,amsthm}
\usepackage{enumitem}
\usepackage{xcolor} % Revisions are marked in blue.
\PassOptionsToPackage{hidelinks,hypertexnames=false}{hyperref}
\usepackage{doi}

\makeatletter
\@namedef{subjclassname@2020}{\textup{2020} Mathematics Subject Classification}
\makeatother

\theoremstyle{plain}
\newtheorem{theorem}{Theorem}[section]
\newtheorem{lemma}[theorem]{Lemma}
\newtheorem{corollary}[theorem]{Corollary}
\newtheorem{proposition}[theorem]{Proposition}

\theoremstyle{definition}
\newtheorem{definition}[theorem]{Definition}

\newtheorem{remark}[theorem]{Remark}
\usepackage{needspace}
\newcommand{\NN}{\mathcal N}
\newcommand{\TT}{\mathcal T}
\newcommand{\DD}{\mathcal D}
\newcommand{\N}{\mathbb N_0}
\newcommand{\V}{\mathsf V}
\newcommand{\C}{\mathcal C}
\newcommand{\supp}{\operatorname{supp}}
\newcommand{\Id}{\operatorname{Id}}
\newcommand{\Sub}{\operatorname{Sub}}
\newcommand{\Mod}{\operatorname{Mod}}
\newcommand{\conv}{\operatorname{conv}}
\newcommand{\cl}{\operatorname{cl}}
\setlist[enumerate]{label=\textup{(\roman*)},leftmargin=*}

\begin{document}
\title[Equational theories of min-plus semirings]
{Equational theories of min-plus semirings}

\author{Yilin Zhou}
\address{School of Mathematics, Northwest University, Xi'an, 710127, Shaanxi, P.R. China}
\email{zhou990619@163.com}

\author{Zidong Gao}
\address{School of Mathematics, Northwest University, Xi'an, 710127, Shaanxi, P.R. China}
\email{zidonggao@yeah.net}

\subjclass[2020]{16Y60, 08B05, 08B15}
\keywords{additively idempotent semiring, limit variety, finite identity basis,
min-plus algebra, subvariety lattice.}

\begin{abstract}
We study the equational theory and subvarieties of the variety
$\mathcal V$ generated by $(\mathbb N_0,\min,+)$.
We prove that a subvariety of $\mathcal V$ is proper if and only if
it is locally finite, equivalently, if and only if it satisfies a
power-stabilization identity. Every finite member generates a
hereditarily finitely based variety with a finite subvariety lattice.
In contrast, the periodic subvarieties defined by $x^m=x^{m+1}$
are nonfinitely based for $m\geq2$.
We construct a strictly descending chain of nonfinitely based
subvarieties and three varieties generated by explicit reduced
products, each with a decidable equational theory. The smallest is
a limit variety: it is nonfinitely based, and every proper
subvariety has an absolute finite identity basis. The subvariety
lattice of $\mathcal V$ has cardinality continuum. In fact, the
interval between the first reduced-product variety and the
subvariety defined by $x^2=x^3$ already has this cardinality,
and all its members are nonfinitely based. We also give a
hereditarily finitely based example with an infinite subvariety
lattice. The constructions use finite convex models and recognition
quotients.
\end{abstract}

\maketitle

\section{Introduction}\label{sec:introduction}
An additively idempotent semiring is an algebra whose addition is
a semilattice operation and whose multiplication is associative
and distributes over addition on both sides. Such algebras are
also called semilattice-ordered semigroups. Examples include
tropical semirings, semirings of relations, and endomorphism
semirings of semilattices; see \cite{Golan1999,Pin1998}.
Throughout the paper we use the signature $(+,\cdot)$ without
constants. Thus an identity or zero element of an algebra, when
present, is not named in the signature.

The finite basis problem and the structure of subvariety lattices
are two basic questions in the study of semiring varieties.
For semigroups, these questions have a long history; see the survey
of Volkov \cite{Volkov2001}. In the semiring setting, the
interaction between the two operations leads to further difficulties. Dolinka constructed a nonfinitely based finite semiring
\cite{Dolinka2007}, developed criteria for inherent nonfinite
basability \cite{Dolinka2009INFB}, and studied related examples
and endomorphism semirings \cite{Dolinka2009Remark,Dolinka2009End}.
Jackson, Ren and Zhao \cite{JRZ2022} subsequently showed that
nonfinite basability can occur even when the multiplicative
semigroup has a finite identity basis.

Several classes of semiring varieties admit explicit structural
descriptions. Polin \cite{Polin1980} studied minimal semiring
varieties. Subvariety lattices or term normal forms have been
described for ordered bands and commutative multiplicatively
idempotent semirings \cite{Pastijn2005,ChajdaLanger2018}.
Other classifications concern varieties generated by two-element
ai-semirings \cite{ShaoRen2015} and commutative ai-semiring
varieties satisfying power identities \cite{RenZhaoShao2020}. The finite basis problem has also been
studied systematically for three-element semirings
\cite{ZhaoEtAl2020}, for several classes of four-element semirings
\cite{RenEtAl2025I,YueEtAl2025II,RenEtAl2026III}, and for extensions
of small nonfinitely based examples \cite{WuRenZhao2024}.
Subvariety lattices can also be large: Dolinka \cite{Dolinka2008}
constructed an interval of semiring varieties of cardinality
continuum.

A \emph{limit variety} is a nonfinitely based variety all of whose
proper subvarieties are finitely based. Ren, Jackson, Zhao and
Lei \cite{RJZL2023} constructed limit varieties of ai-semirings
using flat extensions of groups and a combinatorial example.
They also proved existence results for further limit varieties.
To identify a particular limit variety, however, one must describe
its identities or generators and prove that \emph{every} proper
subvariety is finitely based.

In this paper we study the commutative ai-semiring variety
\[
 A=(\mathbb N_0,\min,+_{\mathrm{ar}}),\qquad
 \mathcal V=\V(A),
\]
where $+_{\mathrm{ar}}$ is ordinary addition of nonnegative integers.
The carrier contains zero but not infinity. Aceto, \'Esik and
Ing\'olfsd\'ottir \cite{AEI2003} gave geometric descriptions of
identities and free algebras for tropical semirings and proved
nonfinite-basis results. We use their convex-geometric description
in our constructions. Metcalfe and Santschi \cite{MS2025} studied
related questions for idempotent semifields, whose multiplication
is invertible; this assumption does not hold for $A$.

Here it is necessary to distinguish min-plus from max-plus.
In semiring notation, $A$ satisfies
\[
 x+xy=x.
\]
The max-plus algebra on the nonnegative integers satisfies the
opposite absorption identity $x+xy=xy$; it is this latter algebra
that occurs in \cite[Section~4.1]{RJZL2023}.
Moreover, if a member of $\mathcal V$ has an element $o$ absorbing
both operations, then $a=a+ao=a+o=o$ for every $a$.
Hence $\mathcal V$ contains no nontrivial flat semiring in the
sense of \cite{RJZL2023}, so a construction in $\mathcal V$
requires different generating algebras.

We first prove that the proper subvarieties of $\mathcal V$ are
exactly those satisfying a power-stabilization identity, and hence
exactly its locally finite subvarieties. Their finite basis
properties nevertheless differ. Every finite member generates a
hereditarily finitely based variety, whereas some locally finite
subvarieties have no basis in a bounded number of variables.
We construct a descending chain of nonfinitely based subvarieties
ending in a limit variety. The finite models used in these
constructions also yield intervals of cardinality continuum.

For $m\geq1$, write
\[
 \mathcal U_m=\mathcal V\cap\Mod(x^m=x^{m+1}).
\]
The finite truncations $A_n,B_n$ are defined in
Section~\ref{sec:B}. Section~\ref{sec:descents} defines successive subvarieties
$\mathcal L_2,\mathcal L_3,
\mathcal L_4$ of $\mathcal U_2$ by explicit identities.
The reduced-product varieties $\mathcal W_\infty$, $\mathcal Y$
and $\mathcal Z$ are constructed in
Sections~\ref{eventual:section}, \ref{smaller:section}
and~\ref{clique:section}, respectively. Finally, set
\[
 \mathcal G=\mathcal U_2\cap\Mod(xy^2=x^2y^2),\qquad
 R_k:\quad x_1\cdots x_k=\sum_{i=1}^k(x_1\cdots x_k)x_i.
\]
The following theorem collects the main results.

\Needspace{7\baselineskip}
\begin{theorem}[Main theorem]\label{limit:main}
In the constant-free signature $(2,2)$, let
$A=(\mathbb N_0,\min,+_{\mathrm{ar}})$ and $\mathcal V=\V(A)$.
\begin{enumerate}
\item A subvariety of $\mathcal V$ is proper if and only if it is locally
finite, equivalently if and only if it satisfies $x^m=x^{m+1}$ for
some $m\geq1$.
\item Every finite member of $\mathcal V$ generates a hereditarily
finitely based variety with a finite subvariety lattice. Some proper subvarieties are contained in no member of either
truncation chain $\V(A_n)$ or $\V(B_n)$, although both chains
have join $\mathcal V$.
\item The varieties $\mathcal U_m$ are nonfinitely based for all $m\geq2$.
There is an explicitly defined strictly descending chain of
nonfinitely based subvarieties
\[
 \mathcal L_4\subsetneq\mathcal L_3\subsetneq\mathcal L_2
 \subsetneq\mathcal U_2\subsetneq\mathcal V.
\]
Neither $\mathcal V$, any $\mathcal U_m$, nor any of
$\mathcal L_2,\mathcal L_3,\mathcal L_4$ is a limit variety.
\item The cofinite reduced product $M$ defined in
Theorem~\ref{eventual:main} generates a locally finite nonfinitely
based variety $\mathcal W_\infty$ with a decidable equational theory.
Moreover $\mathcal W_\infty\subsetneq\mathcal W_*\subseteq\mathcal L_4$,
where $\mathcal W_*$ has the safe-identity presentation in
Section~\ref{separation:section}. There is a further explicitly defined nonfinitely based subvariety
$\mathcal Y\subsetneq\mathcal W_\infty$ with a decidable
equational theory. Thus $\mathcal W_\infty$ is not a limit variety.
The symmetric recognition construction gives
$\mathcal Z\subsetneq\mathcal Y$ with a decidable equational theory.
The variety $\mathcal Z$ is a limit subvariety: it is nonfinitely
based and every proper subvariety is finitely based. Consequently
$\mathcal Y$ is not a limit variety. Every limit subvariety of $\mathcal V$ is locally
finite, lies in some $\mathcal U_m$ with $m\geq2$, and cannot be
generated by a finite algebra or a finite family of finite algebras.
\item The lattices of subvarieties of $\mathcal V$ and of its proper
locally finite subvariety $\mathcal U_2$ each have cardinality
$2^{\aleph_0}$. Each contains continuum many nonfinitely based
subvarieties. The interval $[\mathcal W_\infty,\mathcal U_2]$
already has continuum many distinct nonfinitely based members.
\item The subvariety $\mathcal G$ is hereditarily finitely based and has
an infinite subvariety lattice; it fails every $R_k$. Thus neither an infinite subvariety lattice nor failure of all
the saturation identities $R_k$ guarantees the existence of a
limit subvariety.
\end{enumerate}
We construct one explicit limit subvariety and determine the
cardinality of the full subvariety lattice. The uniqueness and
classification of limit subvarieties remain open here.
\end{theorem}

Our proofs use finite convex models for the periodic varieties.
Saturation identities give finite bases and finite subvariety
lattices for varieties generated by finite algebras. For the
nonfinite-basis results, we construct finite algebras satisfying
all identities in a prescribed number of variables but failing
an identity in more variables. Suitable recognition quotients
then give reduced-product varieties with decidable equational
theories. Strengthening the recognition conditions yields
$\mathcal Z\subsetneq\mathcal Y\subsetneq\mathcal W_\infty$.
We prove that every proper subvariety of $\mathcal Z$ is finitely
based by reducing additional identities to a canonical form.

The paper is organized as follows. Section~\ref{sec:preliminaries}
introduces the convex models. Section~\ref{sec:finite} proves
the finite basis results for finite members and their subvarieties.
Section~\ref{sec:descents} constructs a descending chain of periodic
nonfinitely based varieties. Sections~\ref{eventual:section}
and~\ref{card:section} give the first two reduced-product
constructions and an interval of cardinality continuum.
Section~\ref{clique:section} constructs $\mathcal Z$, whose
limit property is proved in Section~\ref{clmin:section}.
Section~\ref{sec:conclusion} completes the proof of
Theorem~\ref{limit:main}.

\section{Identities and finite convex models}\label{sec:preliminaries}
All algebras are nonempty and have the signature $(+,\cdot)$.
We write multiplication by juxtaposition and use $=$ for an
identity. Let $\C$ be the finite set of axioms for commutative
additively idempotent semirings: associativity, commutativity and
idempotence of addition, associativity and commutativity of
multiplication, and the two distributive laws.
The variety $\mathcal V=\V(A)$ satisfies these axioms and
\[
 I:\quad x+xy=x.
\]
We use $\V(\mathcal K)$ for the variety generated by a class
$\mathcal K$, $\Id(\mathcal K)$ for its identities, and
$\Mod(\Sigma)$ for the class defined by a set of identities
$\Sigma$. We use Birkhoff's variety theorem and standard
equational logic throughout \cite{Birkhoff1935,BurrisSankappanavar1981}.

\subsection{Equational conventions}
\begin{definition}
A variety is \emph{finitely based} if a finite set of identities defines
it among all algebras in this fixed signature. It is \emph{hereditarily
finitely based} if all its subvarieties are finitely based. It is a
\emph{limit variety} if it is nonfinitely based and all its proper
subvarieties are finitely based.
\end{definition}
Thus a finite basis must define the variety among all algebras
in the fixed signature. A finite set of additional identities over
a nonfinitely based ambient variety does not by itself provide
such a basis. A variety is \emph{locally
finite} if all its finitely generated members are finite.

\begin{lemma}[Polynomial representation]\label{found:polynomial}
Every term in $r\geq1$ variables is equivalent modulo $\C$ to a finite
nonempty sum of nonempty monomials. It has an exponent representation
$P\subseteq\N^r\setminus\{\boldsymbol0\}$, finite and nonempty,
with value $\min_{p\in P}p\cdot a$ at $a\in\N^r$ in $A$.
\end{lemma}
\begin{proof}
Represent each variable by its unit vector. Addition corresponds
to union of exponent sets, with repetitions removed by additive
idempotence, and multiplication corresponds to pairwise vector sums
by distributivity. Induction on terms gives the polynomial expansion
and the evaluation formula. Neither operation produces the empty set
or the zero vector.
\end{proof}
Write $t_P=\sum_{p\in P}x^p$ for the polynomial represented by
$P$. The support of a monomial is the set of variables with
positive exponent. A zero exponent means that the corresponding
factor is absent; it does not denote a constant term.
Modulo additive idempotence, two polynomials are equal if and only
if each absorbs every monomial of the other.
We sometimes write $P+x^q=P$ for $t_P+x^q=t_P$.

\begin{lemma}[Exact identity criterion]\label{found:convex}
For a finite nonempty $P\subseteq\N^r\setminus\{0\}$ put
$K(P)=\conv(P)+\mathbb R_{\geq0}^r$. A nonzero exponent vector $q$
satisfies $A\models P+x^q=P$ if and only if $q\in K(P)$.
Consequently $A\models P=Q$ if and only if $K(P)=K(Q)$.
\end{lemma}
\begin{proof}
If $q\geq\sum_{p\in P}\lambda_p p$ with $\lambda_p\geq0$ and
$\sum\lambda_p=1$, then for every $a\in\N^r$,
\[
 q\cdot a\geq\sum_p\lambda_p(p\cdot a)\geq\min_p p\cdot a.
\]
This proves the absorption. Conversely suppose $q\notin K(P)$.
The set $K(P)$ is closed: for a convergent sequence $c_j+u_j$ in it,
compactness of $\conv(P)$ supplies a convergent subsequence of $c_j$;
the corresponding $u_j$ converge in the closed nonnegative cone.
Choose a nearest point $z\in K(P)$ to $q$, which exists by restricting
a minimizing sequence to a closed bounded ball. Put $w=z-q$.
Differentiating squared distance along the segment from $z$ to any
$v\in K(P)$ gives $w\cdot(v-z)\geq0$. Since $z+te_i\in K(P)$
for $t\geq0$, all $w_i\geq0$. Moreover
$w\cdot(z-q)=\|w\|^2>0$, so $w\cdot(p-q)>0$ for all $p\in P$.
The finitely many strict inequalities still hold after a sufficiently
small perturbation to a nonnegative rational vector. Clearing
denominators gives an integer assignment $a\geq0$ such that
$q\cdot a<p\cdot a$ for every $p\in P$, so the absorption fails.
The last assertion follows by applying this criterion to every
monomial on each side.
\end{proof}
This is the constant-free form of the convex-geometric identity
criterion in \cite{AEI2003}. We have included the separation argument
for the carrier and signature used here. For general convex
separation results, see \cite{Rockafellar1970}.

\begin{lemma}[Positive witnesses]\label{found:positive}
The subalgebra $A^+=(\mathbb N_{>0},\min,+)$ has the same identities
as $A$.
\end{lemma}
\begin{proof}
Identities pass to subalgebras. For the converse, suppose $s(a)<t(a)$
at an assignment in $\N^r$, interchanging sides if necessary. Let
$D\geq1$ bound the degrees of all monomials in both terms, choose an
integer $M>D$, and put $b_i=Ma_i+1$. For each exponent $p$,
$p\cdot b=M(p\cdot a)+\deg(p)$. Hence
\[
 s(b)\leq Ms(a)+D<Mt(a)+1\leq t(b),
\]
since the original integer values differ by at least one. Each
$b_i$ is positive, so the identity also fails in $A^+$.
\end{proof}

\begin{theorem}[Proper subvarieties]\label{found:proper}
For a subvariety $\mathcal W\subseteq\mathcal V$, the following are
equivalent: $\mathcal W\ne\mathcal V$; $\mathcal W$ satisfies
$x^m=x^{m+1}$ for some $m\geq1$; and $\mathcal W$ is locally finite.
In the second case its relatively free algebra on $r\geq1$ generators
has size at most $2^{(m+1)^r-1}-1$.
\end{theorem}
\begin{proof}
A proper subvariety has an identity $s=t$ failing in $A$, since otherwise
$A\in\mathcal W$. By Lemma~\ref{found:positive} choose a positive
integer assignment with $s(a)=u<v=t(a)$. Substitute the nonempty term
$x^{a_i}$ for each variable. In $\C+I$, $x^d$ absorbs $x^e$ for
$e\geq d$: use idempotence if $e=d$, and $I$ with the nonempty factor
$x^{e-d}$ if $e>d$. Each unary polynomial therefore reduces to the power with the
least exponent. The substitution yields $x^u=x^v$ with $u\geq1$.
The absorptions $x^u+x^{u+1}=x^u$ and
$x^{u+1}+x^v=x^{u+1}$, together with the equality of the endpoint
powers, imply $x^u=x^{u+1}$.

Under $x^m=x^{m+1}$ every exponent at least $m$ reduces to $m$, by
successive multiplication and induction. Thus an $r$-variable term is
represented by a nonempty subset of the $(m+1)^r-1$ nonzero vectors
in $\{0,\ldots,m\}^r$, proving the bound and local finiteness.
Finally, the subalgebra of $A$ generated by the numerical element $1$
contains every positive integer, obtained by repeated multiplication.
This is an infinite one-generated algebra, so $\mathcal V$ is not
locally finite. This proves the remaining implication.
\end{proof}

\begin{proposition}\label{found:periodiccover}
Put $\mathcal U_m=\mathcal V\cap\Mod(x^m=x^{m+1})$ for $m\geq1$.
These varieties form an increasing chain of proper locally finite
subvarieties, and every proper subvariety of $\mathcal V$ is
contained in some $\mathcal U_m$. Moreover $\mathcal U_1=\mathcal D$ is
the distributive-lattice variety, with precisely the trivial subvariety
and itself; both are finitely based.
\end{proposition}
\begin{proof}
Every proper subvariety lies in one of these varieties by
Theorem~\ref{found:proper}. Multiplying the power identity by $x$
gives the increasing inclusions. For $m=1$, $\C,I,x^2=x$ imply the two lattice absorptions:
$x+xy=x$ and $x(x+y)=x^2+xy=x$. They define distributive lattices.
To check completeness, expand each term as a sum of squarefree
monomials. In the two-element lattice, $P+q=P$ holds exactly when
some monomial of $P$ has support contained in that of $q$: assign the
multiplicative unit to the target support and the other element outside
it for necessity; use $I$ for sufficiency. Thus these finite axioms
generate the two-element lattice variety. That lattice belongs to
$\mathcal V$: send the numerical zero of $A$ to its multiplicative
unit and every positive integer to its other element. Minimum and
ordinary addition show directly that this is a surjective homomorphism. Every nontrivial lattice has a two-element
sublattice: use two distinct comparable elements, or an incomparable
element and its join with another. A nontrivial subvariety consequently
contains the two-element generator. The trivial variety has basis
$x=y$.
\end{proof}

\subsection{Capped convex free algebras}\label{sec:convex}

For positive integers $m,r$, put
\[
 \mathcal U_m=\V(A)\cap\operatorname{Mod}(x^m=x^{m+1}),\qquad
 Q_{m,r}=\{0,\ldots,m\}^r\setminus\{\boldsymbol0\}.
\]
Write $c_m$ for coordinatewise truncation at $m$. For a nonempty
$S\subseteq Q_{m,r}$ define
\[
 \operatorname{cl}(S)=Q_{m,r}\cap
       \bigl(\operatorname{conv}(S)+\mathbb R_{\geq0}^r\bigr).
\]
Let $C_{m,r}$ consist of the nonempty closed subsets, with operations
\[
 S+T=\operatorname{cl}(S\cup T),\qquad
 ST=\operatorname{cl}\{c_m(s+t):s\in S,\ t\in T\}.
\]
Addition of exponent vectors means ordinary vector addition.
Since the zero vector is excluded, each generating vector represents
a nonempty monomial, and no constant term is needed.

\begin{lemma}\label{cap:closure}
For nonempty $S,T\subseteq Q_{m,r}$,
\[
 \operatorname{cl}\{c_m(a+b):a\in\operatorname{cl}(S),\ b\in\operatorname{cl}(T)\}
 =\operatorname{cl}\{c_m(s+t):s\in S,\ t\in T\}.
\]
\end{lemma}
\begin{proof}
If $a\in\operatorname{cl}(S)$, there are coefficients
$\lambda_s\geq0$, with $\sum_s\lambda_s=1$, such that
$a\geq\sum_s\lambda_s s$ coordinatewise. Monotonicity and concavity
of each coordinate of $c_m$ give, for fixed $b$,
\[
 c_m(a+b)\geq c_m\Bigl(\sum_s\lambda_s(s+b)\Bigr)
             \geq\sum_s\lambda_s c_m(s+b).
\]
Thus $c_m(a+b)$ belongs to the convex upward closure of
$\{c_m(s+b):s\in S\}$. Applying the same argument to
$b\in\operatorname{cl}(T)$ proves one inclusion. Extensivity gives
the other. The operator used here is idempotent, since
$\operatorname{conv}(\operatorname{cl}(S))+\mathbb R_{\geq0}^r$
is contained in $\operatorname{conv}(S)+\mathbb R_{\geq0}^r$.
\end{proof}

\begin{proposition}\label{cap:membership}
The algebra $C_{m,r}$ is a finite member of $\mathcal U_m$.
\end{proposition}
\begin{proof}
The set $Q_{m,r}$ is a commutative semigroup under
$(a,b)\mapsto c_m(a+b)$: coordinatewise,
$\min(m,\min(m,a+b)+c)=\min(m,a+b+c)$. Its nonempty powerset is
a commutative ai-semiring under union and setwise product.
Lemma~\ref{cap:closure} says that equality of closures is a congruence
for multiplication; the closure axioms give compatibility with union.
The closed-set quotient is $C_{m,r}$, and it is finite.

To prove membership in $\V(A)$, realize the relatively free algebra
$F_{\mathcal V}(r)$ as an algebra of term functions.
By Lemma~\ref{found:convex}, two exponent sets define the same term
function if and only if their convex upward polyhedra coincide.

Send the term function represented by $P$ to
$\Phi(P)=\operatorname{cl}(c_m(P))$. If
$q\geq\sum_p\lambda_p p$, concavity gives
$c_m(q)\geq\sum_p\lambda_p c_m(p)$. Equality of convex upward
polyhedra consequently implies equality of the images under $\Phi$.
Thus $\Phi$ is well defined on the relatively free algebra
$F_{\V(A)}(r)$, realized as the algebra of term functions in
$A^{A^r}$. It is a homomorphism by Lemma~\ref{cap:closure} and
$c_m(p+q)=c_m(c_m(p)+c_m(q))$. It is surjective, since each closed set is the image of the
polynomial whose monomials have those exponent vectors. Therefore
$C_{m,r}\in\V(A)$.

For every closed $S$ and every positive integer $k$,
\[
 S^k=\operatorname{cl}\{c_m(ks):s\in S\}.
\]
Each sum $c_m(s_1+\cdots+s_k)$ is coordinatewise at least the
average of $c_m(ks_1),\allowbreak\ldots,\allowbreak c_m(ks_k)$, by concavity. This gives one containment; taking $k$ equal vectors gives the
reverse containment.
For an integer vector $s\geq0$,
$c_m(ms)=c_m((m+1)s)$: a coordinate is zero on both sides when the
original coordinate is zero, and is $m$ on both sides otherwise.
The formula therefore gives $S^m=S^{m+1}$.
\end{proof}

\begin{proposition}\label{cap:free}
The elements $X_i=\operatorname{cl}(\{e_i\})$ freely generate
$C_{m,r}$ in $\mathcal U_m$. For exponent sets $P,Q$ of terms in $r$
variables,
\[
 \mathcal U_m\models P=Q
 \quad\Longleftrightarrow\quad
 \operatorname{cl}(c_m(P))=\operatorname{cl}(c_m(Q)).
\]
\end{proposition}
\begin{proof}
The $X_i$ generate every capped monomial and hence every nonempty
closed set. By Proposition~\ref{cap:membership}, every identity of
$\mathcal U_m$ holds at these generators in $C_{m,r}$. This proves
necessity.
Conversely, power stabilization reduces every exponent greater than
$m$ to $m$. Equality of the two displayed closed sets implies that
each capped generating vector belongs to the convex upward polyhedron
of the other capped set. Lemma~\ref{found:convex} gives an identity of $A$ between the capped polynomials. This
identity holds in $\mathcal U_m$, and the exponent reductions prove
the original identity.
\end{proof}

\begin{remark}\label{rem:principal}
For $v\in Q_{m,r}$, write $\mathord\uparrow v=\cl(\{v\})$.
Every element of $C_{m,r}$ is the sum of its principal upper sets.
The power formula in Proposition~\ref{cap:membership} also gives
$(S+T)^k=S^k+T^k$ for every $k\geq1$: apply the concavity of
$c_m$ to each point of $\cl(S\cup T)$ before multiplying by $k$.
In particular, $C_{2,r}$ satisfies $(x+y)^2=x^2+y^2$.
\end{remark}

\section{Finite members and hereditary finite bases}\label{sec:finite}
Not every proper subvariety of $\mathcal V$ is generated by a finite
convex model. Nevertheless, every finite member lies in a finitely
axiomatized variety with a finite subvariety lattice. We prove this
using truncation and saturation identities, and then give a
hereditarily finitely based variety with an infinite subvariety lattice.

\subsection{Truncations and nilpotent identities}\label{sec:B}
For $n\geq1$ set $A_n=\{0,1,\ldots,n,\infty\}$ and
$B_n=\{1,\ldots,n,\infty\}$. Addition is numerical minimum. For finite
$a,b$, multiplication has value $a+b$ if this ordinary integer sum is
at most $n$, and has value $\infty$ otherwise. A product with an
$\infty$ input is $\infty$.
Numerical $0$ in $A_n$ is a multiplicative unit. The element $\infty$
is an additive identity and a multiplicative absorbing element.
The signature still has no constant symbols for these elements.
The semilattice order $a\leq_S b$ defined by $a+b=b$
is the reverse of numerical order.

\begin{lemma}\label{lem:algebras}
Both $A_n$ and $B_n$ satisfy $\C$. They satisfy
\begin{equation}\label{eq:I}x+xy= x.\tag{I}\end{equation}
Every product of $n+1$ elements of $B_n$ is $\infty$.
\end{lemma}
\begin{proof}
Minimum is associative, commutative and idempotent. A finite product has
the ordinary sum of its inputs as value if this is at most $n$, and
$\infty$ otherwise; {an input equal to $\infty$ gives $\infty$}. Since all inputs are nonnegative, truncating an intermediate sum
does not change the final value. The value is therefore independent
of the order and bracketing of the factors. Multiplication by a fixed
element is monotone in the numerical chain, hence preserves the minimum
of two elements. This proves both distributive laws. Multiplication
cannot decrease either input numerically, proving \eqref{eq:I}.
In $B_n$ every finite input is at least $1$, so $n+1$ inputs have sum
greater than $n$.
\end{proof}

Fix $n\geq1$ and put $m=n+1$. Define $\NN_n$ by $\C$ and
\begin{align}
 x_1\cdots x_m&= y^m,\label{eq:N1}\tag{N1}\\
 x+y^m&= x.\label{eq:N2}\tag{N2}
\end{align}

\begin{lemma}[Common zero term]\label{lem:zero}
In $\NN_n$, the value of $z^m$ is independent of $z$, is an additive
identity, and absorbs multiplication. Every product of length at least
$m$ equals this value. Every term is equivalent to a sum of monomials of
degrees $1$ through $n$, or to the term $z^m$.
\end{lemma}
\begin{proof}
Substitute $a$ for all $x_i$ and $b$ for $y$ in \eqref{eq:N1} to obtain
$a^m=b^m$. Write $e$ for this common value in the metatheory.
Equation~\eqref{eq:N2} gives $a+e=a=e+a$. A product of length $L\geq m$
can be grouped into $m$ nonempty factors: the first $m-1$ are single
factors and the last contains the remaining $L-m+1$ factors.
Equation~\eqref{eq:N1} makes its value $e$. In particular $a^m b=e$,
so $eb=be=e$.

Expand a term by distributivity. Each monomial of degree at least $m$
equals $e$ and can be removed if another summand remains. If all monomials have degree at least the bound, the term is
equivalent to $z^m$, where $z$ may be chosen fresh. All reductions
use substitutions and applications of identities in a term context;
no empty product or constant symbol is needed.
\end{proof}

Let $M_n$ be the nonempty monomials in $x_1,\ldots,x_n$ of degree at most
$n$. Their number is
\begin{equation}\label{eq:dn}d_n=\binom{2n}{n}-1.\end{equation}
Indeed, adding a slack coordinate counts weak compositions of $n$ into
$n+1$ parts; the zero exponent vector is then removed. For $P\subseteq M_n$
write $T_P=\sum_{p\in P}p$ when $P\ne\varnothing$ and
$T_\varnothing=z^m$, with $z$ distinct from the $x_i$.
Define the finite pool
\begin{equation}\label{eq:Fpool}
 \mathcal F_n=\{T_P+q= T_P:P\subseteq M_n,\ q\in M_n\}.
\end{equation}
There are $d_n2^{d_n}$ indexed candidates. For $W\subseteq\NN_n$,
let $E_n(W)$ consist of those that hold in the variety.

\begin{lemma}[Completeness of the nilpotent pool]\label{lem:Bcomplete}
For every subvariety $W$ of $\NN_n$,
$\C\cup\{\textup{\eqref{eq:N1}},\textup{\eqref{eq:N2}}\}\cup E_n(W)$
is a finite identity basis for $W$.
\end{lemma}
\begin{proof}
The proposed axioms hold by construction. To prove completeness,
let $s= t$ be any identity of $W$. By Lemma~\ref{lem:zero}, its sides
are represented by sets $P,Q$ of short monomials, with the empty set
representing the common zero term. Choose all auxiliary variables
fresh, distinct from the variables occurring in these monomials.

For $q\in Q$, validity of $s= t$ and idempotence give
$s+q= s$ in $W$. Put $U=\supp(q)$; then $|U|\leq n$.
Replace every original variable outside $U$ by a fresh zero term $z^m$.
Multiplicative absorption deletes precisely those monomials of $P$ whose
support is not contained in $U$. Thus $W$ satisfies
\begin{equation}\label{eq:Brestrict}
 T_{P_U}+q= T_{P_U},\qquad
 P_U=\{p\in P:\supp(p)\subseteq U\}.
\end{equation}
This includes $P_U=\varnothing$. Rename $U$ injectively into
$\{x_1,\ldots,x_n\}$, keeping $z$ fresh. The resulting identity belongs
to $E_n(W)$, so renaming back derives \eqref{eq:Brestrict} from the
proposed finite basis. If $P\setminus P_U$ is nonempty, add its sum to
both sides. The zero-term identity and additive associativity give
$s+q= s$. If there are no omitted monomials, this is already
\eqref{eq:Brestrict}. The argument includes $P=\varnothing$.

Repeating over $Q$ gives $s+t= s$. If $Q$ is empty this follows
directly from the additive identity property. Interchanging the sides
gives $t+s= t$, and hence $s= t$.
\end{proof}

\begin{proposition}\label{prop:Bhereditary}
Every subvariety of $\NN_n$ is finitely based, and
\[
 |\Sub(\NN_n)|\leq 2^{d_n2^{d_n}}.
\]
In particular $\V(B_n)$ is hereditarily finitely based and has a finite
subvariety lattice. A complete finite basis for $B_n$ is effectively
obtained by retaining the true members of \eqref{eq:Fpool} and adjoining
$\C$, \eqref{eq:N1} and \eqref{eq:N2}.
\end{proposition}
\begin{proof}
Lemma~\ref{lem:Bcomplete} gives each finite basis. It also makes
$W\mapsto E_n(W)$ injective, since the same subset gives the same
axiomatization. There are at most $2^{d_n2^{d_n}}$ such subsets.
Lemma~\ref{lem:algebras} shows $B_n\in\NN_n$; minimum with $\infty$
also verifies \eqref{eq:N2}. Hence every subvariety of $\V(B_n)$ is
covered. To decide a pool identity, test the $(n+1)^n$ assignments to
$x_1,\ldots,x_n$; the value of $z^m$ is always $\infty$.
This finite test determines which identities to retain, and
Lemma~\ref{lem:Bcomplete} shows that the resulting basis is complete.
\end{proof}

\subsection{Saturation and a variety-theoretic dichotomy}\label{sec:T}
For $n\geq1$ and $m=n+1$, define $\TT_n$ by $\C$, \eqref{eq:I}, and
\begin{align}
 x^m&= x^{m+1},\label{eq:P}\tag{P}\\
 x_1\cdots x_m&=\sum_{i=1}^m x_i(x_1\cdots x_m).
 \label{eq:R}\tag{R}
\end{align}
Let $D$ be the algebra on $\{0,\infty\}$ with minimum as addition and
maximum as multiplication, and put $\DD=\V(D)$. This is the two-element
distributive lattice in our operation order. Direct inspection shows
$D\in\TT_n$: its multiplication is idempotent, the lattice absorption
law gives \eqref{eq:I}, and every summand in \eqref{eq:R} equals its
left side.

\begin{lemma}[Saturation normal form]\label{lem:saturation}
Modulo the axioms of $\TT_n$, every term is a finite nonempty sum of
nonempty monomials with positive exponents at most $m$, such that the
sum of the exponents strictly below $m$ is at most $n$.
The normalization preserves the support of each expanded monomial.
\end{lemma}
\begin{proof}
Expand by distributivity. Multiplying \eqref{eq:P} by additional copies
of its variable reduces every exponent at least $m$ to $m$.
In a monomial of support $U$, call an exponent equal to $m$ saturated.
If the sum of nonsaturated exponents is at least $m$, choose $m$
occurrences among those coordinates, allowing repeated variables.
Use these occurrences as the $m$ inputs of \eqref{eq:R}, then multiply
by the remaining factors if any. Each resulting summand increases one
selected nonsaturated exponent by one. No exponent exceeds $m$, and
the support remains $U$.

On every branch the nonnegative integer
$\sum_{u\in U}(m-\operatorname{exp}_u)$ strictly decreases.
Strong induction on this integer yields a finite expansion with no
further applicable step. Each terminal monomial has nonsaturated degree
at most $m-1=n$. Thus the procedure terminates and gives a normal form equal to
the original term modulo the axioms.
\end{proof}

\begin{lemma}[Divisibility and support restriction]\label{lem:Arestrict}
In $\TT_n$, a monomial absorbs every multiple of itself.
Suppose $W\subseteq\TT_n$, $D\in W$, and $P+q= P$ holds in $W$,
where $P$ is a nonempty polynomial and $q$ is a nonempty monomial with
all exponents at most $m$. Exponents of $P$ may also be capped at $m$.
For $U=\supp(q)$, the set $P_U$ of summands supported in $U$ is
nonempty, and $P_U+q= P_U$ follows by substitution and the
axioms of $\TT_n$ from the original absorption.
\end{lemma}
\begin{proof}
If $r=pv$ with $v$ nonempty, \eqref{eq:I} gives $p+r= p$.
If $r=p$, the same equality follows from idempotence. Thus one
polynomial absorbs another whenever each summand of the latter is a
multiple of a summand of the former.

In $D$, assign $0$ to $U$ and $\infty$ outside $U$. The target $q$ has
value $0$, whereas $P$ has value $0$ precisely when $P_U$ is nonempty.
Thus validity forces nonemptiness. Set
$e_U=\prod_{u\in U}u^m$, an actual term because $U$ is nonempty.
Substitute $e_U$ for every variable outside $U$. A summand in $P_U$
is unchanged. Any other summand acquires all variables of $U$ to
exponent at least $m$, and reduces to $e_U$ by \eqref{eq:P}.
Each member of $P_U$ divides $e_U$ and hence absorbs it.
The substituted polynomial therefore reduces to $P_U$, while $q$
remains unchanged. This does not require $e_U$ to be a common zero.
\end{proof}

Let $\mathcal M_n$ be the nonempty monomials on $n+1$ fixed variables with
exponents in $\{0,\ldots,m\}$. Put
\[
 h_n=|\mathcal M_n|=(n+2)^{n+1}-1,
\]
\[
 \mathcal G_n=\left\{\sum_{p\in P}p+q=\sum_{p\in P}p:
 \varnothing\ne P\subseteq \mathcal M_n,\ q\in \mathcal M_n\right\}.
\]
There are $h_n(2^{h_n}-1)$ candidate identities. Let $G_n(W)$ consist
of those that hold in $W$. Both sides of each candidate are nonempty
polynomials.

\begin{lemma}[Completeness above $\DD$]\label{lem:Acomplete}
If $D\in W\subseteq\TT_n$, the axioms of $\TT_n$ together with
$G_n(W)$ form a finite identity basis for $W$.
\end{lemma}
\begin{proof}
All the proposed axioms hold in $W$. Let $s= t$ be an identity of
$W$ and normalize both sides by Lemma~\ref{lem:saturation}. For each monomial
$q$ of $t$, validity gives $s+q= s$.
Lemma~\ref{lem:Arestrict} gives a nonempty $P_U$ and a valid absorption
$P_U+q= P_U$, where $U=\supp(q)$.
Split $U=X\mathbin{\dot\cup}Y$, with the nonsaturated exponents on $X$
and exponent $m$ on $Y$. Then $|X|\leq n$.

If $Y$ is empty, the restricted absorption has at most $n$ variables
and capped exponents. After an injective renaming, it belongs to
$G_n(W)$. Reversing the renaming and adding the omitted summands of
$s$ gives $s+q= s$.

If $Y$ is nonempty, first identify all its variables with a fresh $z$
in the valid restricted absorption, multiply both sides by $z^m$,
and cap exponents. The resulting absorption has at most $n+1$
variables, capped exponents and nonempty sides, so is in $G_n(W)$
after renaming. Its target is the $X$-part of $q$ times $z^m$, or
just $z^m$ if $X$ is empty. From this selected axiom, rename $X$
back and substitute $\prod_{y\in Y}y$ for $z$. Writing
$g=\prod_{y\in Y}y^m$, we obtain
\[
 H+q= H,\qquad H=\sum_{p\in P_U}pg,
\]
modulo exponent reduction. The factor $z^m$ ensures the presence of
every saturated variable, including in summands originally omitting it.
Every $pg$ is a multiple of a summand of $s$, so
$s+H= s$ by \eqref{eq:I}. Add $s$ to the displayed absorption
to derive $s+q= s$. The argument also applies when $X=\varnothing$, without introducing
a constant term. Deriving all the target absorptions and then
interchanging the two sides proves $s= t$.
\end{proof}

\begin{lemma}[The alternative excluding $D$]\label{lem:dichotomy}
If $W\subseteq\TT_n$ and $D\notin W$, then $W\subseteq\NN_n$.
Within $\TT_n$, exclusion of $D$ is equivalent to \eqref{eq:N2}.
\end{lemma}
\begin{proof}
Choose an identity $s= t$ valid in $W$ and failing in $D$,
and normalize its sides. Interchanging them if needed, choose an
assignment with $s=\infty$ and $t=0$. Some monomial $q$ of $t$
has value $0$. For $U=\supp(q)$, every summand of $s$ has a variable
outside $U$, since otherwise that summand also has value $0$.
The absorption $s+q= s$ is valid in $W$.

Substitute $y$ for all variables in $U$ and a distinct $x$ for all
others. The target becomes $y^d$ for some $d\geq1$, and every
summand on the left is divisible by $x$. Adding $x$ and using
\eqref{eq:I} gives $x+y^d= x$. If $d\geq m$, cap its exponent.
If $d<m$, divisibility gives $y^d+y^m= y^d$, whence
\[
 x+y^m=(x+y^d)+y^m= x+(y^d+y^m)= x.
\]
Thus \eqref{eq:N2} holds. The values $a^m$ all coincide, because
$a^m+b^m=a^m=b^m$. Denote the common additive identity by $e$.
For any $b$, $eb=b^m b=b^{m+1}=b^m=e$; commutativity gives $be=e$.

Normalize $x_1\cdots x_m$ with distinct variables by
Lemma~\ref{lem:saturation}. Every terminal summand still has support
of size $m$. It must have a saturated variable, since otherwise its
nonsaturated degree is at least $m$. That variable contributes $e$,
so each summand equals $e$ and their sum equals $e$. This proves
\eqref{eq:N1}. Conversely, \eqref{eq:N2} fails in $D$ at
$x=\infty,y=0$, proving the stated equivalence.
\end{proof}

\begin{proposition}[Ambient hereditary theorem]\label{prop:T}
The variety $\TT_n$ is hereditarily finitely based, with
\begin{equation}\label{eq:latticebound}
 |\Sub(\TT_n)|\leq
 2^{h_n(2^{h_n}-1)}+2^{d_n2^{d_n}}.
\end{equation}
Each subvariety has a finite basis using at most $n+2$ variables
in each identity.
\end{proposition}
\begin{proof}
Subvarieties containing $D$ have the bases in Lemma~\ref{lem:Acomplete}
and are determined by their subsets of $\mathcal G_n$. Their number
is bounded by the first summand. Those excluding $D$ lie in $\NN_n$
by Lemma~\ref{lem:dichotomy} and satisfy the second bound by
Proposition~\ref{prop:Bhereditary}. This proves both finiteness assertions. The first pool uses $n+1$
variables, and the corresponding defining axioms use at most
$\max(3,n+1)$. The nilpotent pool uses at most $n+1$ variables,
including its fresh zero-term variable; \eqref{eq:N1} uses $n+2$.
The axioms $\C$ use at most three. Since $n\geq1$, the stated bound
covers all cases.
\end{proof}

\subsection{Membership of the finite truncations}\label{sec:A}
\begin{lemma}\label{lem:Amember}
For every $n\geq1$, $A_n\in\TT_n$.
\end{lemma}
\begin{proof}
Lemma~\ref{lem:algebras} verifies $\C$ and \eqref{eq:I}.
For \eqref{eq:P}, a numerical zero input gives value $0$ on both
sides, whereas a positive input or $\infty$ gives $\infty$.
{For \eqref{eq:R}, an input equal to $\infty$ makes both sides equal to $\infty$.}
If all $m=n+1$ inputs are positive finite integers, their sum exceeds
$n$, so both sides again equal $\infty$. Otherwise one input is
numerical $0$. Multiplication by that input leaves the original
product unchanged, while the other summands are numerically no
smaller. Their minimum is therefore the original product.
\end{proof}

\subsection{Finite members and the truncation chains}

\begin{corollary}\label{cap:counterexample}
For every $m\geq2$ and $r\geq1$, $C_{m,r}$ fails $R_r$.
Consequently $\mathcal U_m$ is contained in no single $\TT_n$,
$\V(A_n)$, or $\V(B_n)$.
\end{corollary}
\begin{proof}
At the principal generators, the left side of $R_r$ is the closure of
$\delta=(1,\ldots,1)$ and the right side is the closure of
$\{\delta+e_i:1\leq i\leq r\}$. Every point in the latter convex
upward polyhedron has coordinate sum at least $r+1$, while $\delta$
has sum $r$. Thus the values differ. Each of the three stated finite
stages satisfies $R_{n+1}$, giving the exclusions.
\end{proof}

\begin{proposition}\label{cap:finite-saturation}
The algebra $C_{m,r}$ satisfies $R_{(m-1)r+1}$.
\end{proposition}
\begin{proof}
Put $k=(m-1)r+1$ and evaluate at closed sets $S_1,\ldots,S_k$.
A generating vector of their product is
$a=c_m(a_1+\cdots+a_k)$, where $a_i\in S_i$. If the $j$th coordinate
of $a$ is less than $m$, at most $m-1$ of the factors have a positive
$j$th coordinate. Thus at most $(m-1)r$ factors have a positive entry in at least one
unsaturated coordinate. Some $a_i$ is therefore zero in
every unsaturated coordinate of $a$. Repeating this factor changes
no coordinate: $c_m(a+a_i)=a$. Hence $a$ occurs among the generators
of the $i$th right-hand summand of $R_k$. Closure gives containment
of the left side in the right side. Conversely, each right-hand
generator is coordinatewise above a left-hand generator, so upward
closure gives the reverse containment.
\end{proof}

\begin{theorem}\label{cap:finite-members}
Every finite algebra $S\in\V(A)$ generates a hereditarily finitely
based variety with a finite subvariety lattice. If $|S|=s\geq2$, then
\[
 \V(S)\subseteq\TT_{(s-1)s}.
\]
No limit subvariety of $\V(A)$ is generated by a finite algebra or
by a finite family of finite algebras.
\end{theorem}
\begin{proof}
The identity $x+xy=x$ makes the power sequence of each element
descending in the semilattice order $a\leq b\iff a+b=b$.
A descending sequence in a set of $s$ elements stabilizes by the
$s$th power; equality of consecutive powers propagates under
multiplication. Thus $S\models x^s=x^{s+1}$ and $S\in\mathcal U_s$.
Taking all its elements as generators, we obtain from
Proposition~\ref{cap:free} a surjective homomorphism $C_{s,s}\to S$. Proposition~\ref{cap:finite-saturation} shows that
$S$ satisfies $R_{(s-1)s+1}$. With $N=(s-1)s$, the power stabilization
at $s$ also gives stabilization at $N+1$. These identities, together
with $\C$ and $x+xy=x$, put $S$ in $\TT_N$. Both conclusions now follow from the hereditary finite-basis and
finite-lattice theorem for $\TT_N$. For $s=1$, the generated variety is trivial.

A finite family of finite algebras generates the same variety as
their finite direct product: the product belongs to the variety of
the family, and every factor is a surjective homomorphic image of
the product. The first assertion applies to that product. Since a limit variety is nonfinitely based, neither type of finite
generating family is possible.
\end{proof}

\begin{proposition}\label{limit:joins}
The chains $\V(A_n)$ and $\V(B_n)$ are strictly increasing, and
\[
 \mathcal V=\bigvee_{n\geq1}\V(A_n)=\bigvee_{n\geq1}\V(B_n).
\]
Neither set-theoretic union is a variety. Each finite stage is
hereditarily finitely based and has a finite subvariety lattice.
\end{proposition}
\begin{proof}
The map $\pi_n:A\to A_n$ fixes $0,\ldots,n$ and sends larger
integers to $\infty$. It preserves minimum by monotonicity. A sum
$a+b\leq n$ forces both summands to be at most $n$ and is unchanged;
a sum greater than $n$ is sent to $\infty$, as is the product of
the capped inputs. Thus $\pi_n$ is a surjective homomorphism.
Its restriction maps $A^+$ onto $B_n$, and $B_n$ is also a
subalgebra of $A_n$. Truncation from stage $n+1$ to $n$ gives a
surjective homomorphism in each family. The identity $x^{n+1}=x^{n+2}$ holds at stage $n$ and
fails at the element $1$ at stage $n+1$, proving strictness.

Every identity failing in $A$ has an assignment with distinct finite
term values. Taking $n$ at least their maximum preserves their
difference under $\pi_n$, since term evaluation commutes with
homomorphisms by induction. Thus the common identities of all the
$A_n$ are exactly $\Id(A)$. Using a positive witness from Lemma~\ref{found:positive} gives the
same conclusion for the $B_n$.
This proves both join equalities. However $A$ itself satisfies no
identity $x^{n+1}=x^{n+2}$, as the assignment $x=1$ shows. It belongs
to neither union. A union that were a variety would contain the join
of its members, which is $\mathcal V$ and contains $A$, a contradiction.
The final assertion follows from Lemma~\ref{lem:Amember} and
Propositions~\ref{prop:Bhereditary} and~\ref{prop:T}.
\end{proof}
\subsection{A hereditarily finitely based auxiliary variety}\label{guarded:section}

Let $\mathsf B_1$ consist of the finite axioms $\C$ and
\[
 I:\ x+xy=x,\qquad P:\ x^2=x^3,\qquad
 G:\ xy^2=x^2y^2.
\]
Put $\mathcal G=\mathcal U_2\cap\operatorname{Mod}(G)$.

\begin{lemma}\label{guarded:normal}
Modulo $\mathsf B_1$, each monomial is squarefree or is the square
of a squarefree monomial with the same support. If $q=x_1\cdots x_k$,
$k\geq2$, and $i\neq j$, then
\[
 q+\Bigl(\prod_{h\neq i}x_h\Bigr)^2
  +\Bigl(\prod_{h\neq j}x_h\Bigr)^2
 =\Bigl(\prod_{h\neq i}x_h\Bigr)^2
  +\Bigl(\prod_{h\neq j}x_h\Bigr)^2.
\]
\end{lemma}
\begin{proof}
$P$ caps positive exponents at two. A monomial with a repeated variable
has the form $y^2t$. If $t$ is nonempty, $G$ turns it into $y^2t^2$;
after capping, this is the square of the support product. If there is no remaining factor, the monomial is already $y^2$.

Applying $P$ to $x+y$, expanding the square and cube, absorbing
$x^2y$ and $xy^2$ by $I$, and capping $x^3,y^3$, gives
\[
 x^2+y^2+xy=x^2+y^2. \tag{F}
\]
For $k=2$ this is the assertion. For $k>2$, multiply (F) by the
nonempty product of the other $k-2$ variables and apply $G$ to the
two nonlinear summands. This yields the displayed identity.
\end{proof}

\begin{theorem}\label{guarded:basis}
The finite set $\mathsf B_1$ is an absolute identity basis for
$\mathcal G$ in the constant-free signature.
\end{theorem}
\begin{proof}
Every member of $\mathcal G$ satisfies $\mathsf B_1$. To prove the
reverse inclusion, it suffices to derive every identity of $A$ from
$\mathsf B_1$. Expand such an identity into polynomials and consider
the monomial absorptions in both directions. Cap the exponents at two, as permitted by $P$.
The resulting absorption $P_0+q=P_0$ remains valid in $A$: the convex
criterion and concavity of the cap map preserve membership of $q$
in the relevant convex upward polyhedron.

Let $U$ be the support of $q$, and let $P_U$ be the monomials of $P_0$
supported in $U$. This set is nonempty: assign numerical zero on $U$
and numerical one elsewhere. Further, $P_U+q=P_U$ is valid in $A$.
For any assignment on $U$, give every outside variable an integer
larger than the value of $q$. Every outside monomial then has value greater than that of $q$, so some
member of $P_U$ must have value at most that of the target.

If $q$ has a repeated variable, $G$ turns it into the all-squared
support monomial $e_U$. Every capped member of $P_U$ divides $e_U$,
and $I$ proves the absorption. If $q$ is squarefree and $P_U$ contains
a squarefree monomial, that monomial divides $q$. It remains to
treat squarefree $q$ when every member of $P_U$ is nonlinear.
Consider
\[
 M_U=\bigcup_{p\in P_U}(U\setminus\operatorname{supp}(p)).
\]
If $M_U$ is empty, the all-one assignment makes every $p$ strictly
larger than $q$. If $M_U=\{i\}$, assign one to $x_i$ and two to
the other variables of $U$. A monomial with full support has a
doubled coordinate and is strictly larger than $q$. A monomial
missing $i$ has a doubled coordinate elsewhere, whose increase of
at least two exceeds the omitted contribution of one. Both cases
contradict validity. Hence $M_U$ contains distinct $i,j$.

Choose members omitting $i,j$, respectively. After $G$, their squared
support monomials divide $e_{U\setminus\{i\}}$ and
$e_{U\setminus\{j\}}$. Thus $P_U$ absorbs these two terms by $I$,
and Lemma~\ref{guarded:normal} implies that it absorbs $q$.
Adding back discarded summands gives $P_0+q=P_0$. Repeating this in both directions derives the original identity
from $\mathsf B_1$. Hence
$\operatorname{Mod}(\mathsf B_1)\subseteq\V(A)$, and the identities
$P,G$ give equality with $\mathcal G$.
\end{proof}

Recall that $\mathcal D$ is the distributive-lattice variety. Define
\[
 J_k:\ x_1\cdots x_k=(x_1\cdots x_k)^2\quad(k\geq1),
\]
and
\[
 F_k:\ x_1\cdots x_k+(x_2\cdots x_k)^2
       =(x_2\cdots x_k)^2\quad(k\geq2).
\]

\begin{lemma}\label{guarded:canonical}
If $\mathcal D\subseteq\mathcal W\subseteq\mathcal G$, then
$\mathcal W$ is axiomatized by $\mathsf B_1$ together with the
identities $J_k,F_k$ that hold in $\mathcal W$. Within $\mathsf B_1$,
$J_k$ implies $J_l$ for $l\geq k$, and $F_k$ implies $F_l$ for $l\geq k$.
\end{lemma}
\begin{proof}
Normalize the monomials of an identity of $\mathcal W$ using
Lemma~\ref{guarded:normal}, and consider an absorption $P_0+q=P_0$.
For $U=\operatorname{supp}(q)$, the subset $P_U$ of summands supported
in $U$ is nonempty. In the two-element lattice, give variables in $U$
the top value and outside variables the bottom value; validity forces
a summand supported in $U$.

The absorption $P_U+q=P_U$ also holds in $\mathcal W$. Replace all
outside variables by the nonempty term
$e_U=(\prod_{i\in U}x_i)^2$. Every discarded summand becomes $e_U$
under $G,P$, and every member of $P_U$ divides $e_U$ and absorbs it
under $I$. The substituted identity is therefore the desired one.
If no terms were discarded, no substitution is required.

An all-squared $q$, or a squarefree $q$ with a squarefree divisor in
$P_U$, gives an absorption already following from $I$. Otherwise
$q$ is squarefree and all members of $P_U$ are squared support
monomials. An omitted-coordinate union of size at least two is
handled by Lemma~\ref{guarded:normal}. If that union is empty,
$P_U=e_U$, and the absorption is $J_{|U|}$ because also $q+e_U=q$.
If it is $\{i\}$, $P_U$ reduces by $I$ to $e_{U\setminus\{i\}}$,
and the absorption is a renaming of $F_{|U|}$. The monomial with smaller support is nonempty, so $|U|\geq2$.
This exhausts the cases. Adding the discarded summands back and
applying the argument in both directions proves the claimed
axiomatization.

For $l>k$, multiply $J_k$ by the nonempty product of the extra
variables. $G$ turns the right side into the square of the full
product, giving $J_l$. The same multiplication and $G$ turn $F_k$
into $F_l$. For $l=k$ there is nothing to prove.
\end{proof}

\begin{theorem}\label{guarded:hfb}
The variety $\mathcal G$ is hereditarily finitely based. Every
subvariety containing $\mathcal D$ has a basis consisting of
$\mathsf B_1$ and at most two additional identities.
\end{theorem}
\begin{proof}
For a subvariety containing $\mathcal D$, select the least index of
a valid $J_k$, if one exists, and the least index of a valid $F_k$,
if one exists. The implication chains and Lemma~\ref{guarded:canonical} show that
these choices give a finite basis. It is absolute because
$\mathsf B_1$ is finite.

Suppose $\mathcal W$ excludes $\mathcal D$. An identity of
$\mathcal W$ failing in the two-element lattice yields an absorption
$P_0+q=P_0$ and $U=\operatorname{supp}(q)$ such that every summand of
$P_0$ contains an outside variable. Send all variables of $U$ to $y$
and all others to $x$. The target is $y^d$ for some $d\geq1$.
Adding $x$ and using $I$ gives $x+y^d=x$. For $d=1$ commutativity
makes the variety trivial. For $d\geq2$, $P$ gives $x+y^2=x$.
Thus all squares have a common value $e$, which is an additive
identity. Identity (F) gives $e+ab=e$, and therefore $ab=e$ for all
$a,b$.

Every term in such an algebra is equivalent either to a nonempty
sum of variables or to the common product value. In a nontrivial algebra choose $a\neq e$.
Then $\{a,e\}$ is a two-element subalgebra: addition has
identity $e$, and every product equals $e$. Denote this algebra
by $E$. Conversely a missing linear variable
is detected by assigning it $a$ and all other variables $e$;
a nonempty linear sum is distinguished from the common product in
the same way. These observations show that every nontrivial algebra
of this kind has exactly the identities of $E$. Hence
$\mathcal W$ is trivial or $\V(E)$, both finitely based. The trivial variety has basis $x=y$;
the variety $\V(E)$ is defined by $\mathsf B_1$ together with
$x+y^2=x$. The preceding term description proves completeness
of this latter finite basis.
\end{proof}

\subsection{A support model and an infinite lattice}

Fix $r\geq1$. For each nonempty $S\subseteq\{1,\ldots,r\}$ take
atoms $\ell_S,e_S$. Set
\[
 e_Se_T=e_{S\cup T},\qquad e_S\ell_T=e_{S\cup T},\qquad
 \ell_S\ell_T=
 \begin{cases}\ell_{S\cup T},&S\cap T=\varnothing,\\
 e_{S\cup T},&S\cap T\neq\varnothing.
 \end{cases}
\]
For a nonempty set $P_0$ of atoms define $d(P_0)$ by the following rules:
$e_U\in d(P_0)$ exactly when a member of $P_0$ has support contained
in $U$; and $\ell_U\in d(P_0)$ exactly when either $P_0$ contains
$\ell_S$ with $S\subseteq U$, or
\[
 \left|\bigcup_{e_S\in P_0,\ S\subseteq U}(U\setminus S)\right|\geq2.
\]
If the indexing family is empty, this set-theoretic union is empty;
no empty algebraic sum is involved.

\begin{proposition}\label{guarded:free-support}
The nonempty $d$-closed atom sets, with closure of union and closure
of setwise product, form the free algebra of $\mathcal G$ on $r$
generators.
\end{proposition}
\begin{proof}
A product of atoms is linear exactly when all factors are linear
and have pairwise disjoint supports; otherwise it is squared on
the union of supports. This proves associativity and commutativity.

The closure $d$ is extensive and monotone. For idempotence, replace
each newly generated squared atom on support $S$ by an original atom
on a subset of $S$. If such an original atom is linear, it already
forces every relevant larger linear support. Otherwise replacement
by the smaller squared support can only enlarge the omitted-coordinate
union. A newly generated linear atom has either an original linear
subset-support or two omitted coordinates supplied by original
squared atoms. Enlargement of its support preserves these conditions.
Applying the closure a second time therefore adds no new atom.

For compatibility with multiplication, multiply a newly generated
atom on $U$ by an atom on $V$. If the target product is squared,
an original atom supported in $U$ produces an atom supported in
$U\cup V$, which forces the squared target. If the target product
is linear, the factors are linear and $U,V$ are disjoint. An
original linear subset-support then produces a linear subset-support
of $U\cup V$. In the other case the squared atoms supplying two
missing coordinates in $U$ continue to supply those same missing
coordinates after adjoining $V$. This proves compatibility; the
reverse containment follows from extensivity. Thus closure equality
is a semiring congruence.

$I$ follows from support inclusion and the product rule. The square
of a polynomial is the closure of the squared atoms on its supports.
A mixed product which is squared is absorbed by one squared factor
support. A linear mixed product has two disjoint nonempty supports,
whose squared atoms supply at least two missing coordinates and
therefore absorb it. The same description gives $P_0^3=P_0^2$:
multiplication produces squared unions, and repetition of the same
support recovers every original squared support. It also gives
$P_0Q_0^2=P_0^2Q_0^2$, since both are generated by squared unions
of one support from each factor. The model satisfies $\mathsf B_1$.

The principal singleton-support linear atoms generate every atom
and closed set. Every absorption specified by $d$ follows from $\mathsf B_1$:
the squared case and the first linear case follow from divisibility,
and the second linear case follows from Lemma~\ref{guarded:normal}
and divisibility. Any absorption not specified by $d$ fails at these
generators in the model. Thus $d$ determines exactly the identities
of $\mathsf B_1$, so freeness follows from
Theorem~\ref{guarded:basis}.
\end{proof}

\begin{corollary}\label{guarded:infinite-lattice}
The hereditarily finitely based variety $\mathcal G$ fails every $R_k$ and has an infinite
subvariety lattice. More precisely,
$\mathcal G\cap\operatorname{Mod}(J_k)$, $k\geq1$, is a strictly
increasing chain.
\end{corollary}
\begin{proof}
In the support model on $k$ generators the full product is $\ell_U$,
where $U$ is the full support, whereas its square is $e_U$.
The closure of $\{e_U\}$ omits $\ell_U$, since its support misses
no coordinate. Hence $J_k$ fails. In $\mathsf B_1$, each summand
on the right side of $R_k$ equals the square of the product: for
$k\geq2$ apply $G$ to its repeated variable, and for $k=1$ this
is the displayed identity itself. Thus $R_k$ and $J_k$ are equivalent
in $\mathcal G$.

In the model on $k$ generators, a product of $k+1$ atoms cannot
have pairwise disjoint nonempty supports, and so is squared.
A product of $k+1$ arbitrary elements is therefore a sum of squared
atoms and is idempotent by the description of squares in
Proposition~\ref{guarded:free-support}. This model satisfies
$J_{k+1}$ and fails $J_k$. Together with the implication
$J_k\Rightarrow J_{k+1}$, this proves every strict inclusion.
\end{proof}

\section{Nonfinite bases and proper descents}\label{sec:descents}
\subsection{Punctured convex models}
\label{nfb:section}

Put $\mathcal U_2=\mathcal V\cap\operatorname{Mod}(x^2=x^3)$.
We first construct finite algebras that distinguish identities according to
the number of variables they involve.  The construction modifies the punctured
convex-filter argument of Aceto, \'Esik and Ing\'olfsd\'ottir
\cite[Lemma~35, Proposition~45, and the proof of Theorem~32,
pp.~441--447]{AEI2003}.  The indecomposability statement used in their proof no longer holds
when coordinates are capped.  We replace it with the closure argument
below.

Throughout this section $k\geq3$ and
\[
 Q_k=\{0,1,2\}^k,\qquad
 c(v)_i=\min\{2,v_i\},\qquad
 C(S)=Q_k\cap\bigl(\conv(S)+\mathbb R_{\geq0}^k\bigr)
\]
for nonempty $S\subseteq Q_k$.  Vector addition in these definitions is
ordinary coordinatewise addition.  The zero vector is an element of the
finite model; it is not a constant symbol of the language.  Indices on
coordinates are read modulo $k$.  Define
\[
 \delta=(1,\ldots,1),\qquad
 \gamma_i=\delta+e_i-e_{i+1},\qquad
 \Gamma=C(\{\gamma_1,\ldots,\gamma_k\}),\qquad
 \Delta=\Gamma\setminus\{\delta\}.
\]
The weight of a vector is the sum of its coordinates.

\begin{lemma}\label{nfb:geometry}
Every vector of $\Gamma$ has weight at least $k$.  Its weight-$k$
vectors are precisely $\delta,\gamma_1,\ldots,\gamma_k$, and the
$\gamma_i$ are extreme points of that face.  If $S\subseteq\Gamma$ and
$\delta\notin S$, then
\[
 \delta\in C(S)
 \quad\Longleftrightarrow\quad
 \{\gamma_1,\ldots,\gamma_k\}\subseteq S.
\]
\end{lemma}

\begin{proof}
Every generator of $\Gamma$ has weight $k$.  A convex combination has
the same weight, and a nonzero nonnegative upward displacement strictly
increases it.  Thus a weight-$k$ vector belongs to the convex hull of
the $\gamma_i$.  Every coordinate of a point in this hull lies between zero and two.
If an integer vector has coordinate two at position $i$, only
$\gamma_i$ can occur in its convex representation, since every other
generator has coordinate less than two there.  Otherwise all its coordinates are at
most one, and weight $k$ forces it to equal $\delta$.  The same
coordinate argument proves that each $\gamma_i$ is extreme.

In a representation
$\delta=\sum_i\lambda_i\gamma_i$, comparison of coordinate $i$ gives
$\lambda_i=\lambda_{i-1}$.  The coefficients sum to one, so every
coefficient is $1/k$.  If $\delta$ is absent, generating it therefore requires all the $\gamma_i$.
Conversely, their average is
$\delta$.
\end{proof}

\begin{lemma}\label{nfb:closure}
Let $D_k$ consist of all nonempty $C$-closed subsets of $Q_k$, together
with $\Delta$.  {It is closed under intersections of nonempty families.}  Its closure
on nonempty subsets is
\[
 C_*(S)=
 \begin{cases}
  \Delta,&C(S)=\Gamma\text{ and }\delta\notin S,\\
  C(S),&\text{otherwise}.
 \end{cases}
\]
Every member of $D_k$ is an upper set in the coordinatewise order.
\end{lemma}

\begin{proof}
No vector of $\Gamma$ other than $\delta$ lies below $\delta$, by
Lemma~\ref{nfb:geometry}.  Removing $\delta$ therefore preserves upward
closure.  Intersections of convex upper sets are convex upper sets.
For a $C$-closed set $F$,
\[
 F\cap\Delta=(F\cap\Gamma)\setminus\{\delta\}.
\]
If $F\cap\Gamma=\Gamma$, this intersection is $\Delta$.  If the
intersection is a proper subset of $\Gamma$, removing $\delta$
preserves $C$-closedness: regenerating it would require all the
$\gamma_i$, and hence all of $\Gamma$, by
Lemma~\ref{nfb:geometry}.  This proves the intersection assertion.

The additional set $\Delta$ can be a closed superset of $S$ only when
$S\subseteq\Delta$, in which case $C(S)\subseteq\Gamma$.
Lemma~\ref{nfb:geometry} shows that a proper $C$-closed subset of
$\Gamma$ cannot newly generate $\delta$.  Thus the displayed exceptional case is the only one in which the new
closure differs from the ordinary closure. The formula therefore gives
the least member of $D_k$ containing $S$.
\end{proof}

For nonempty subsets put
\[
 S\odot T=\{c(s+t):s\in S,\ t\in T\}.
\]

\begin{lemma}\label{nfb:compatibility}
For all nonempty $S,T\subseteq Q_k$,
\[
 C_*\bigl(C_*(S)\odot C_*(T)\bigr)=C_*(S\odot T).
\]
\end{lemma}

\begin{proof}
First, ordinary closure $C$ is compatible with $\odot$.  If
$a\geq\sum_s\lambda_s s$ coordinatewise, monotonicity and concavity of
coordinatewise capping give
\[
 c(a+b)\geq c\Bigl(\sum_s\lambda_s(s+b)\Bigr)
          \geq\sum_s\lambda_s c(s+b).
\]
Applying this argument to both factors proves compatibility for $C$.
Moreover, $C(C_*(S))=C(S)$.  The two raw products in the lemma therefore
have the same ordinary closure; denote it by $F$.  If
$F\ne\Gamma$, both $C_*$-closures equal $F$.  Suppose $F=\Gamma$.
By the exceptional closure rule, we need only show that if
$\delta$ belongs to the larger raw product, then it also belongs to
$S\odot T$.

Suppose $a\in C_*(S)$, $b\in C_*(T)$, and $c(a+b)=\delta$.
Then $a+b=\delta$, since a coordinate sum at least two cannot cap to
one.  Ordinary compatibility implies
\[
 S\odot\{b\}\subseteq\Gamma,
 \qquad \delta\in C(S\odot\{b\}).
\]
If $b\ne0$, choose $j$ with $b_j>0$.  Every vector in
$S\odot\{b\}$ is at least $b$, so this set does not contain
$\gamma_{j-1}$, whose $j$th coordinate is zero.
Lemma~\ref{nfb:geometry} forces $\delta\in S\odot\{b\}$.
Choose $s\in S$ with $c(s+b)=\delta$.  If $s\ne0$, the same
fixed-translate argument, now applied to $\{s\}\odot T$, gives
$t\in T$ with $c(s+t)=\delta$.

If instead $s=0$, then $b=\delta$ and $0\in S$.  Since
$S\odot T\subseteq\Gamma$, it follows that $T\subseteq\Gamma$.
The facts $\delta\in C_*(T)$ and $C(T)\subseteq\Gamma$ force
$\delta\in T$: a proper closed subset cannot newly generate it by
Lemma~\ref{nfb:geometry}, and the exceptional closure $\Delta$ omits
it.  Thus $\delta\in S\odot T$ in this case as well.

Finally, if $b=0$, then $0\in T$.  Indeed a nonnegative convex
combination with a nonnegative upward displacement can equal zero only
if zero was already in the original set.  Now $S\subseteq\Gamma$ and
$\delta=a\in C_*(S)$ force $\delta\in S$ by the preceding reasoning.
Again $\delta\in S\odot T$.  The reverse inclusion of raw products follows from extensivity.
This proves the equality.
\end{proof}

\begin{proposition}\label{nfb:witnesses}
Define operations on $D_k$ by
\[
 S+T=C_*(S\cup T),\qquad ST=C_*(S\odot T).
\]
Then $D_k$ is a finite commutative additively idempotent semiring
satisfying $x+xy=x$ and $x^2=x^3$.  Ordinary closure gives a
homomorphism
\[
 h:D_k\longrightarrow \widehat C_k,\qquad
 h(\Delta)=\Gamma,\quad h(S)=S\ (S\ne\Delta),
\]
where $\widehat C_k$ is the algebra of nonempty ordinary convex upper sets in
$Q_k$.  The constant-free algebra $\widehat C_k$ belongs to $\mathcal U_2$.
The algebra $D_k$ fails the identity
\begin{equation}\label{nfb:cyclic}
 E_k:\quad r_k+s_k=s_k,\qquad
 r_k=\prod_{j=1}^k x_j,\qquad
 s_k=\sum_{i=1}^k x_i^2
                  \prod_{j\notin\{i,i+1\}}x_j.
\end{equation}
This identity holds in $A$.
\end{proposition}

\begin{proof}
Capped vector addition makes $Q_k$ a commutative semigroup.  Its
nonempty powerset under union and complex multiplication is a
commutative additively idempotent semiring.
Lemma~\ref{nfb:compatibility} says that closure equivalence for $C_*$ is
a multiplicative congruence; the closure axioms give compatibility with
union.  The quotient by closure equivalence is $D_k$, so the semiring laws
hold.  It is finite because $Q_k$ is finite.  All its
elements are upper sets, and $c(s+t)\geq s$, so $ST\subseteq S$ and
$S+ST=S$.  Ordinary closure produces $\widehat C_k$ and the displayed
homomorphism $h$.

The presence of the zero exponent vector in $\widehat C_k$ does not
enlarge the signature.  To verify membership in the constant-free
variety, consider the subalgebra of $A^{A^k}$ consisting of the functions
\[
 a\longmapsto\min_{p\in P}p\cdot a
 \quad(P\subseteq{\N^k}\text{ finite and nonempty}),
\]
where the zero exponent vector is permitted.  The constant-zero
function is an element of this direct power.  The ordinary convex
criterion and the concavity calculation in
Lemma~\ref{nfb:compatibility} show that
$P\mapsto C(c(P))$ is a well-defined surjective homomorphism onto
$\widehat C_k$.  Thus its constant-free reduct belongs to $\mathcal V$.

For every ordinary closed $S$ and integer $d\geq1$,
\begin{equation}\label{nfb:ordinarypower}
 S^d=C(\{c(ds):s\in S\}).
\end{equation}
By concavity, a capped sum of $d$ vectors is at least the average
of their $d$-fold capped multiples.  Conversely, each $d$-fold
multiple is obtained by taking $d$ equal vectors.  Taking closures
proves \eqref{nfb:ordinarypower}.  For integer vectors,
$c(2s)=c(3s)$, giving square stabilization in $\widehat C_k$.

An ordinary square $S^2$ cannot equal $\Gamma$.  Its generators in
\eqref{nfb:ordinarypower} have all coordinates in $\{0,2\}$.  If their
closure were $\Gamma$, every extreme $\gamma_i$ would occur among
them, by Lemma~\ref{nfb:geometry}.  But $\gamma_i$ has a coordinate
one because $k\geq3$.  For $T\in D_k$, the common value
$h(T)^2=h(T)^3$ therefore differs from $\Gamma$.  Every other value
has a unique $h$-preimage, so $T^2=T^3$.

At $X_i=C_*(\{e_i\})$, the polynomial $s_k$ evaluates to $\Delta$,
whereas $r_k+s_k$ evaluates to $\Gamma$.  This proves failure of
\eqref{nfb:cyclic} in $D_k$.  In $A$ the identity follows from
$\delta=k^{-1}\sum_i\gamma_i$ and the convex criterion.  All products
in \eqref{nfb:cyclic} are nonempty because $k\geq3$.
\end{proof}

\begin{lemma}\label{nfb:smallvariables}
Every constant-free identity of $\mathcal U_2$ on fewer than $k$
variables holds in $D_k$.
\end{lemma}

\begin{proof}
Polynomial expansion and exponent capping at two are valid in $D_k$ by
Proposition~\ref{nfb:witnesses}.  The free-algebra criterion of
Proposition~\ref{cap:free} reduces the assertion to an absorption
$P+q=P$, where
\[
 \varnothing\ne P\subseteq\{0,1,2\}^r\setminus\{0\},\qquad
 q\in\{0,1,2\}^r\setminus\{0\},\qquad r<k,
\]
and $q\in\conv(P)+\mathbb R_{\geq0}^r$.  Choose
$\lambda_p\geq0$ with
\begin{equation}\label{nfb:coefficients}
 \sum_p\lambda_p=1,\qquad \sum_p\lambda_pp\leq q.
\end{equation}

Suppose the absorption fails at $S_1,\ldots,S_r\in D_k$.  Applying
$h$ gives equality in $\widehat C_k\in\mathcal U_2$.  Its only distinct
preimages are $\Gamma$ and $\Delta$.  Since addition is closure of
union, the value of $P$ must be $\Delta$ and that of $P+q$ must be
$\Gamma$.  Consequently every raw vector of $P$ lies in $\Delta$,
every raw vector of $q$ lies in $\Gamma$, and some raw vector of $q$
equals $\delta$.  For the last assertion, the exceptional closure rule shows that a
closure contained in $\Gamma$ can contain $\delta$ only if it was
already in the raw set.

Select $u_{j,l}\in S_j$, $1\leq l\leq q_j$, that give this raw
vector.  Then
\begin{equation}\label{nfb:partition}
 \sum_{j,l}u_{j,l}=\delta.
\end{equation}
No coordinate sum is capped in this equality, since a sum at least
two cannot cap to one.  Thus the selected vectors are zero-one
vectors with pairwise disjoint supports covering all $k$ coordinates.

If $S_j$ contains zero, upper closure implies $S_j=Q_k$.  Every
selected $u_{j,l}$ for this variable is zero.  Otherwise replacing all
of its selections by zero would give a raw vector of $q$ strictly below
$\delta$, of weight less than $k$, contrary to its membership in
$\Gamma$.  Ignore the empty supports of these unit-valued variables.
Each selected vector for any other active variable has nonempty support.

For $p$ with $\lambda_p>0$, equation \eqref{nfb:coefficients} implies
$p_j=0$ whenever $q_j=0$.  In each of the $p_j$ occurrences of
variable $j$, independently choose one of its $q_j$ selected vectors
with equal probability.  For unit-valued variables choose zero.  Every
sample is a raw vector of $P$ and lies in $\Delta$.  Each coordinate
belongs to one selected support and $p_j\leq2$, so the sampled sums
require no capping.  Their average for fixed $p$ is
\[
 \sum_{j:q_j>0}\frac{p_j}{q_j}
                      \sum_{l=1}^{q_j}u_{j,l}.
\]
Average further with the weights $\lambda_p$.
Equations \eqref{nfb:coefficients} and \eqref{nfb:partition} put the
result at most $\delta$ coordinatewise.  Every sample has weight at
least $k$.  Their average must therefore equal $\delta$, and every
positively weighted sample has weight $k$.
{Lemma~\ref{nfb:geometry} and exclusion of $\delta$ imply that every
sample is a $\gamma_i$.  Their weighted average is $\delta$, whose
unique convex representation in these vectors has every coefficient
equal to $1/k$.  Hence all the $\gamma_i$ occur.}

Coordinates in the same selected nonempty support have equal values
in every sample.  If such a support contained distinct coordinates $a,b$,
the sample $\gamma_a$ would have value two at $a$ and a different value
at $b$, a contradiction.  Hence each nonempty selected support is a
singleton.

Suppose a nonunit-valued variable $j$ has $q_j=2$.  Its selected
vectors are $e_a,e_b$ with $a\ne b$.  Choose any $p$ with
$\lambda_p>0$ and fix choices for all other variables.  If $p_j=0$,
the sample has zero at both $a,b$, whereas a $\gamma_i$ has exactly
one zero.  If $p_j=1$, two choices give pairs $(1,0)$ and $(0,1)$ at
$a,b$, with all other coordinates fixed.  Their unique coordinate two
would occur at the same position outside $\{a,b\}$.  A $\gamma_i$
is determined by that position, so the two different samples cannot
both be such vectors.  If $p_j=2$, choices give pairs $(2,0)$ and
$(0,2)$.  They could both be cyclic $\gamma$-vectors only if
$b=a+1$ and $a=b+1$ modulo $k$, which is impossible for $k\geq3$.
These exhaust the possibilities for $p_j$.

Thus every active nonunit-valued variable has $q_j=1$ and contributes
one singleton support.  Covering $k$ coordinates therefore requires
$k$ distinct variables, contrary to $r<k$.  This proves the absorption,
and applying this result to each monomial on both sides proves the
identity.
\end{proof}

\begin{theorem}\label{nfb:U2}
The variety $\mathcal U_2$ is a proper locally finite nonfinitely based
subvariety of $\mathcal V$.  Its identities have no basis with a
bounded number of variables.
\end{theorem}

\begin{proof}
The identity $x^2=x^3$ fails in $A$ at the natural number one, so the
inclusion is proper.  In $r$ variables every term is a nonempty
polynomial in the finitely many nonzero vectors of $\{0,1,2\}^r$.
The free algebra on $r$ generators is finite, proving local finiteness.
If identities on at most $r$ variables formed a basis, choose
$k\geq3$ with $k>r$.  Lemma~\ref{nfb:smallvariables} says $D_k$
satisfies every such identity of $\mathcal U_2$, while
Proposition~\ref{nfb:witnesses} says it fails the identity $E_k$ of
$\mathcal U_2$.  This contradicts soundness of equational deduction.  Since the number
of variables in a finite set of identities is bounded, no finite basis
exists.
\end{proof}

\begin{corollary}\label{nfb:upperinterval}
Every variety $\mathcal W$ satisfying
$\mathcal U_2\subseteq\mathcal W\subseteq\mathcal V$ has no
bounded-variable identity basis.  In particular $\mathcal U_m$ is
nonfinitely based for every $m\geq2$, and $\mathcal V$ is not a limit
variety.
\end{corollary}

\begin{proof}
Every identity of $\mathcal W$ holds in $\mathcal U_2$, and each
$E_k$ holds in $\mathcal V$ and hence in $\mathcal W$.  The same
$D_k$ therefore give the contradiction in the proof of
Theorem~\ref{nfb:U2}.  Square stabilization implies stabilization at
every later exponent, so $\mathcal U_2\subseteq\mathcal U_m$.
Finally, $\mathcal U_2$ is a proper nonfinitely based subvariety of
$\mathcal V$, so $\mathcal V$ is not a limit variety.
\end{proof}

\subsection{A smaller guarded subvariety}

Define
\begin{equation}\label{nfb:L2definition}
 \mathcal L_2=\mathcal U_2\cap\operatorname{Mod}(G_2),
 \qquad G_2:\ z^2xy=z^2x^2y+z^2xy^2.
\end{equation}
We now take quotients of the preceding algebras.  The following lemma
will ensure that the two values separated by the puncture remain
distinct in these quotients.

\begin{lemma}\label{nfb:recognition}
A vector $v\in Q_k$ belongs to $\Gamma$ exactly when either all its
coordinates are positive, or precisely one coordinate $j$ is zero and
the predecessor coordinate $j-1$ is two.  For nonzero $x,y,z\in Q_k$
and arbitrary $w\in Q_k$, put
\[
 a=c(x+y+2z+w),\qquad b=c(a+x),\qquad d=c(a+y).
\]
Then $a\in\Delta$ if and only if $b,d\in\Delta$.
\end{lemma}

\begin{proof}
A vector with no zero lies above $\delta$.  If its $j$th coordinate is
zero, a convex combination of the $\gamma_i$ below it can use only
$\gamma_{j-1}$, the unique generator with zero there.  Membership is
therefore equivalent to lying above that generator.  This proves the
description, including exclusion of a second zero coordinate.

The nonzero vector $z$ ensures that $a$ has a coordinate two, so
$a\ne\delta$.  Also $b,d\geq a$ and $b,d\ne\delta$.
If $a\in\Delta$, upward closure gives $b,d\in\Delta$.
If $a\notin\Gamma$, its zero coordinates remain zero in both $b,d$:
all summands defining $a$ were zero there.  With at least two zeros,
neither $b$ nor $d$ belongs to $\Gamma$.  In the remaining case there
is one zero $j$, and the predecessor coordinate of $a$ is one.  At that
coordinate the uncapped sum $x+y+2z+w$ equals one.  At most one of
$x,y$ is positive there, so at least one of $b,d$ retains value one
and remains outside $\Gamma$.  This proves the reverse implication.
\end{proof}

\begin{lemma}\label{nfb:quotient}
Let $D_k^-=D_k\setminus\{Q_k\}$, and let $\theta_k$ be its least
congruence containing every evaluation instance of $G_2$.
Then $D_k^-$ is a subalgebra, and $\Gamma,\Delta$ have distinct
$\theta_k$-classes.
\end{lemma}

\begin{proof}
Neither addition nor multiplication of nonunit elements can introduce
the zero vector, so removing $Q_k$ gives a subalgebra.  Every element
is a finite sum of principal upper sets of nonzero vectors.  The
Frobenius identity $(x+y)^2=x^2+y^2$ holds in $D_k$ by
Lemma~\ref{nfb:smallvariables}: it is a two-variable identity of
$\mathcal U_2$, and $k\geq3$.  Squares of finite sums consequently
expand diagonally.

For $S\in D_k^-$ let $\sigma(S)$ denote the predicate
$S\subseteq\Delta$.  Since $\Delta$ is $C_*$-closed,
\[
 \sigma(S+T)\quad\Longleftrightarrow\quad
 \sigma(S)\ \text{and}\ \sigma(T).
\]
Declare $S\equiv T$ when $\sigma(SR)$ and $\sigma(TR)$ agree for
every $R\in D_k^-$ and also in the identity context, meaning $S,T$
without an extra multiplication.  This is an equivalence relation.
Associativity proves compatibility with multiplication, since products
of parameters remain in $D_k^-$.  Distributivity and the displayed
conjunction rule give compatibility with addition:
\[
 \sigma((S+U)R)\quad\Longleftrightarrow\quad
 \sigma(SR)\ \text{and}\ \sigma(UR).
\]
The same argument applies in the identity context.  Thus $\equiv$ is
a congruence, and it distinguishes $\Gamma,\Delta$ because their
$\sigma$-values differ in that context.

Expand the three elements in an evaluation of $G_2$, and a
multiplicative parameter, as finite sums of nonzero principal vectors
$x,y,z,w$.  Frobenius gives diagonal expansions of the squares.  For
each choice of these vectors, the three monomials give the vectors
$a,b,d$ of Lemma~\ref{nfb:recognition}.  In the identity context take
the parameter vector to be zero as notation for absence of that
factor; no language constant is adjoined.  Membership of every
left-hand raw vector in $\Delta$ is therefore equivalent to membership
of every right-hand raw vector in $\Delta$.  The closure of a raw
polynomial set lies in $\Delta$ exactly when every raw vector does,
because $\Delta$ is closed.  Thus every instance of $G_2$ belongs to
$\equiv$.  It follows that $\theta_k\subseteq\equiv$, proving the
required separation.
\end{proof}

\begin{theorem}\label{nfb:L2}
The explicitly defined variety $\mathcal L_2$ in
\eqref{nfb:L2definition} is a proper locally finite nonfinitely based
subvariety of $\mathcal U_2$.  It has no bounded-variable identity
basis.  Consequently $\mathcal U_2$ is not a limit variety.
\end{theorem}

\begin{proof}
For $k\geq3$, consider the finite quotient
$\widehat H_k=D_k^-/\theta_k$.  By construction it satisfies $G_2$.
It fails the cyclic identity $E_k$, since $\Gamma,\Delta$ remain
distinct by Lemma~\ref{nfb:quotient}.

Every subalgebra of $D_k^-$ generated by fewer than $k$ elements belongs
to $\mathcal U_2$.  Indeed, express the assigned elements of an arbitrary
identity as terms in the generators.  Substitution gives an identity of
$\mathcal U_2$ on fewer than $k$ variables, which holds by
Lemma~\ref{nfb:smallvariables}.  Every subalgebra of $\widehat H_k$ generated
by fewer than $k$ elements is a homomorphic image of such a subalgebra:
choose preimages of its generators and take the subalgebra they
generate.  It belongs to $\mathcal U_2$ and also satisfies $G_2$,
being a subalgebra of $\widehat H_k$.  It consequently belongs to $\mathcal L_2$.

Each assignment of fewer than $k$ variables in $\widehat H_k$ lies in one of
these generated subalgebras.  Thus $\widehat H_k$ satisfies every such identity
of $\mathcal L_2$, but fails $E_k$, an identity of $\mathcal V$ and
hence of $\mathcal L_2$.  Soundness of equational consequence excludes
every bounded-variable basis and therefore every finite basis.

The inclusion in $\mathcal U_2$ is proper.  At the three free
generators of its ordinary capped convex model, the two sides of $G_2$
are generated respectively by
\[
 (1,1,2)\qquad\text{and}\qquad (2,1,2),\ (1,2,2).
\]
The former has weight four, while every vector in the convex upward
hull of the latter has weight at least five.  Thus $G_2$ fails there.
Local finiteness is inherited from $\mathcal U_2$.  The existence of its
proper nonfinitely based subvariety $\mathcal L_2$ shows that
$\mathcal U_2$ is not a limit variety.
\end{proof}

\subsection{Two further strict descents}
\label{descent:section}

The nonfinitely based variety $\mathcal L_2$ has two further proper
nonfinitely based subvarieties, which we construct below.  In particular,
nonfinite basability alone does not establish the minimality required
of a limit variety.

In five variables, set
\[
 q=xyzwv,\qquad p_1=x^2y^2z^2,\qquad p_2=xw^2v^2,
 \qquad
 \alpha:\ q+p_1+p_2=p_1+p_2,
\]
and put $\mathcal L_3=\mathcal L_2\cap\operatorname{Mod}(\alpha)$.
Let $U=\{x,y,z,w,v\}$ and $D=\{x,y,z\}$, and define
\[
 e_D=x^2y^2z^2,\qquad
 H_i=x_i\prod_{j\in U\setminus\{i\}}x_j^2\quad(i\in D).
\]
Here $x_i$ denotes the variable named by $i$.  Finally, set
\[
 \beta:\ q+e_D+\sum_{i\in D}H_i=e_D+\sum_{i\in D}H_i,
 \qquad
 \mathcal L_4=\mathcal L_3\cap\operatorname{Mod}(\beta).
\]
All sums and products displayed here are nonempty terms.

\begin{lemma}\label{descent:recognition}
Consider an ordinary capped convex algebra or a punctured algebra
$D_k^-$, and let $K$ be one of its closed subsets.  Write
$\sigma_K(S)$ for $S\subseteq K$.  Define $S\equiv_K T$ when
$\sigma_K(SR)$ and $\sigma_K(TR)$ agree for every multiplicative
parameter $R$, and also in the identity context without such a
parameter.  Then $\equiv_K$ is a semiring congruence.

For an absorption $P+q=P$ with exponents at most two, the following
condition puts every instance in $\equiv_K$: for every tuple of
nonzero principal vectors and every exterior multiplicative context,
membership in $K$ of all the raw vectors of $P$ implies membership in
$K$ of the raw vector of $q$.
\end{lemma}

\begin{proof}
Extensivity and closedness of $K$ imply that the closure of a raw set
is contained in $K$ exactly when every raw point belongs to $K$.
In particular,
\[
 \sigma_K(S+T)\quad\Longleftrightarrow\quad
 \sigma_K(S)\ \text{and}\ \sigma_K(T).
\]
For multiplicative compatibility of $\equiv_K$, associativity allows
us to replace the testing parameter by its product with the additional
element.  Additive compatibility follows from distributivity and the
displayed conjunction rule.  Both arguments also apply to the identity
context.  Hence the relation is a congruence.

Every element is a finite sum of principal upper sets.  Frobenius
$(u+v)^2=u^2+v^2$ holds in the ordinary capped algebras and in $D_k$
by Lemma~\ref{nfb:smallvariables}, because $k\geq3$.  Thus arbitrary assignments can be expanded into principal-vector
evaluations, with each square expanded diagonally.  Suppose the whole evaluated
polynomial $P$, with an exterior multiplicative factor, is contained
in $K$.  Then every common tuple of principal choices gives all its
right-hand raw vectors in $K$.  The stated implication puts every raw
$q$ vector in $K$, so adding $q$ preserves the predicate.  If the
evaluated $P$ is not contained in $K$, extensivity ensures that adding
$q$ does not make it contained.  For a variable missing from a monomial of $P$, choose any point from
its nonempty assigned set; this does not alter that monomial or the
order of the universal quantifiers.  An absent
exterior factor is represented by a zero vector only within this
calculation, without introducing a constant term.
\end{proof}

\begin{lemma}\label{descent:safe}
Let $q=x_1\cdots x_r$, and let $P$ be a nonempty set of nonempty
monomials on these variables, with every exponent at most two.
Suppose that the intersection of the sets of doubled variables of
the monomials of $P$ is empty, and that at least one of these
monomials omits at least two variables.  Then every instance of
$P+q=P$ belongs to $\equiv_\Delta$ on every $D_k^-$.
In particular, this conclusion holds for both $\alpha$ and $\beta$.
\end{lemma}

\begin{proof}
Fix a common principal assignment by nonzero vectors and an exterior
context vector, and let $a$ be the capped target vector.  We show that
if $a\notin\Delta$, then at least one right-hand vector lies outside
$\Delta$.

If $a$ has at least two zero coordinates, all summand vectors and the
context vanish there.  Every right-hand vector retains both zeros
and lies outside $\Delta$.  If $a$ has one zero at position $j$ and
predecessor value one, that value comes either from the context alone
or from one unique assigned variable.  In the context case every
right-hand vector retains predecessor value one.  In the variable
case the empty doubled-variable intersection supplies a monomial
whose exponent at that variable is at most one.  Its vector retains
the zero at $j$ and has predecessor value at most one.  It lies
outside $\Gamma$, hence outside $\Delta$.

By Lemma~\ref{nfb:recognition}, the last bad case is $a=\delta$.
All assigned nonzero vectors then have pairwise disjoint nonempty
supports, with entries one, and the context has disjoint support.
A monomial omitting two variables has zeros on two disjoint
nonempty supports and therefore lies outside $\Delta$.
This proves the pointwise implication.  By
Lemma~\ref{descent:recognition}, it holds for arbitrary assignments
and contexts.

For $\alpha$, the doubled-variable sets are $\{x,y,z\}$ and
$\{w,v\}$.  For $\beta$, they are $D$ and $U\setminus\{i\}$ for
$i\in D$, whose intersection is empty.  In both identities the
monomial $e_D$ omits $w,v$.  Thus the two hypotheses apply.
\end{proof}

\begin{lemma}\label{descent:alpha_strict}
The identity $\alpha$ fails in a finite member of $\mathcal L_2$.
Consequently $\mathcal L_3\subsetneq\mathcal L_2$.
\end{lemma}

\begin{proof}
Work in the ordinary free algebra $C_{2,5}$ of
Proposition~\ref{cap:free}, with coordinate order $x,y,z,w,v$.  For the
following vectors write
\[
 p_1=(2,2,2,0,0),\qquad
 p_2=(1,0,0,2,2),\qquad
 \delta=(1,1,1,1,1).
\]
Upward closures in this proof are relative to
$\{0,1,2\}^5\setminus\{0\}$.  Define
\[
 K_1=\mathord{\uparrow}p_1\ \cup\
     \mathord{\uparrow}p_2\ \cup\
     \bigl(\mathord{\uparrow}\delta\setminus\{\delta\}\bigr).
\]
We first show that $K_1$ is convex upward closed in the integer cube.
Let $a$ be an integer vector lying above a convex combination of
points of $K_1$.  A zero coordinate of $a$ cannot be $x$, which is
positive throughout $K_1$.  A zero at $y$ or $z$ forces every point
with positive coefficient into $\mathord{\uparrow}p_2$; a zero at
$w$ or $v$ forces every such point into
$\mathord{\uparrow}p_1$.  The convex combination and the vector
above it then lie in that principal upset.  If $a$ has no zero,
the only excluded integer point is $\delta$.  Every point of $K_1$
has weight at least five, and its unique weight-five point is $p_2$.
Thus a convex combination below $\delta$ would have to equal $p_2$,
which is not below $\delta$.  This excludes the remaining point and
proves closedness.

For a nonlinear vector $a$, meaning a vector with some coordinate two,
put $T=\operatorname{supp}(a)$ and $J(a)=\{i:a_i=1\}$.  Set
\[
 S_1=D_1=\{x,y,z\},\qquad
 S_2=\{x,w,v\},\qquad D_2=\{w,v\},
 \qquad I(T)=\{i:S_i\subseteq T\}.
\]
Its membership in $K_1$ is equivalent to
\begin{equation}\label{descent:K1criterion}
 I(T)\ne\varnothing,
 \qquad
 J(a)\subseteq\bigcup_{i\in I(T)}(T\setminus D_i).
\end{equation}
If only $S_1$ is contained, this requires precisely the coordinates
of $p_1$ to have their prescribed values.  If only $S_2$ is contained,
it requires precisely $a\geq p_2$.  If both are contained, $T=U$
and the two complements cover $U$, so every nonlinear positive
vector qualifies.  If neither is contained, none of the three sets defining $K_1$
contains $a$.  This proves
\eqref{descent:K1criterion}.

Consider a principal instance of $G_2$ with an exterior context:
\[
 a=c(u+v+2t+w),\qquad b=c(a+u),\qquad d=c(a+v),
 \qquad t\ne0.
\]
All three vectors are nonlinear, have the same support $T$, and satisfy
\begin{equation}\label{descent:residualunion}
 J(b)\cup J(d)=J(a).
\end{equation}
Indeed, at a coordinate of $a$ equal to one, $u$ and $v$ cannot both
contribute positively, so repeating at least one retains value one.
A zero or two in $a$ cannot become a residual coordinate.
The fixed-support criterion \eqref{descent:K1criterion} and
\eqref{descent:residualunion} give
$a\in K_1$ exactly when $b,d\in K_1$.
Lemma~\ref{descent:recognition} proves that all $G_2$ instances lie
in $\equiv_{K_1}$.

The finite recognition quotient of $C_{2,5}$ therefore belongs to
$\mathcal U_2$ and satisfies $G_2$, and hence belongs to $\mathcal L_2$.
Its identity context distinguishes the values of $p_1+p_2$ and
$q+p_1+p_2$: the former is contained in $K_1$, whereas the latter
contains $\delta$.  Thus $\alpha$ fails in this quotient.
\end{proof}

\begin{lemma}\label{descent:beta_strict}
The identity $\beta$ fails in a finite member of $\mathcal L_3$.
Consequently $\mathcal L_4\subsetneq\mathcal L_3$.
\end{lemma}

\begin{proof}
In the same integer cube let
\[
 K_0=\mathord{\uparrow}(2,2,2,0,0)
     \ \cup\
     \bigl(\mathord{\uparrow}\delta\setminus\{\delta\}\bigr).
\]
Every point has weight at least six, excluding $\delta$ of weight
five.  If an integer vector above a convex combination has a zero,
every point with positive coefficient must belong to the principal
upset.  Their average and the vector above it still belong to that
upset.  Every other possible integer vector is positive and nonlinear
and already lies in $K_0$.  This proves convex upward closedness.

For a nonlinear vector of support $T$ and residual set $J$, membership
in $K_0$ means either $T=U$, or
$D\subseteq T$ and $J\cap D=\varnothing$.
In a principal $G_2$ instance, support is fixed and the residual sets
satisfy \eqref{descent:residualunion}.  Thus membership of the left
vector is equivalent to membership of both right vectors.
Lemma~\ref{descent:recognition} proves invariance under all $G_2$
instances and contexts.

We check $\alpha$ for the same recognition congruence.  In a common
principal assignment and exterior context, denote its target vector
by $a$ and its two right-hand vectors by $b_1,b_2$.
Suppose $a\notin K_0$.  If $a$ has a zero and omits a coordinate of
$D$, both $b_1,b_2$ retain that omitted coordinate and lie outside
$K_0$.  If $a$ has a zero, contains $D$, and has value one at a
coordinate of $D$, that one is contributed by the context alone or
by one unique variable.  In the context case both right-hand vectors
retain it.  Otherwise the disjoint doubled-variable sets of $p_1,p_2$
leave one right-hand vector with value at most one at that coordinate.
It also retains the zero of $a$, so it is outside $K_0$.

The remaining case is $a=\delta$.  The five assigned nonzero vectors
then have disjoint nonempty supports, and the context is disjoint
from them.  The vector $b_1$ omits the blocks assigned to $w,v$;
$b_2$ omits those assigned to $y,z$.  Each right-hand vector thus
has zero coordinates.  For both to belong to $K_0$, both would have
to lie in $\mathord{\uparrow}(2,2,2,0,0)$.  All four pairwise
disjoint nonempty omitted blocks would then have to lie outside the
three fixed coordinates $D$.  There are only two such coordinates,
which is impossible.  At least one right-hand vector is outside
$K_0$.

By this pointwise implication and
Lemma~\ref{descent:recognition}, every instance of $\alpha$ belongs
to $\equiv_{K_0}$.  Thus the recognition quotient of $C_{2,5}$
belongs to $\mathcal L_3$.  It refutes $\beta$: its
right-hand exponent vectors $e_D,H_x,H_y,H_z$ all lie in $K_0$,
but its left side additionally contains $\delta$.  The identity
context distinguishes the two quotient values.
\end{proof}

\begin{theorem}\label{descent:chain}
There are strict inclusions
\[
 \mathcal L_4\subsetneq\mathcal L_3\subsetneq\mathcal L_2.
\]
Both $\mathcal L_3$ and $\mathcal L_4$ are locally finite and have
no bounded-variable identity basis.  In particular they are
nonfinitely based, and neither $\mathcal L_2$ nor $\mathcal L_3$
is a limit variety.
\end{theorem}

\begin{proof}
For each $k\geq3$, quotient $D_k^-$ by the congruence generated by
all instances of $G_2,\alpha$, and, for the second construction,
also $\beta$.  Lemma~\ref{nfb:quotient} and
Lemma~\ref{descent:safe} put all these relations inside
$\equiv_\Delta$.  Therefore $\Gamma$ and $\Delta$ remain
different, and each quotient fails the cyclic identity $E_k$.

Choose fewer than $k$ generators in either quotient and take preimages
in $D_k^-$.  By Lemma~\ref{nfb:smallvariables}, these preimages
generate a subalgebra in $\mathcal U_2$: writing assigned elements as
terms in the chosen generators reduces an arbitrary identity test to
one on fewer than $k$ variables.  Its homomorphic image also belongs to
$\mathcal U_2$ and satisfies all the additional defining identities,
being a subalgebra of the respective quotient.  It therefore belongs
to $\mathcal L_3$ or $\mathcal L_4$.

Every assignment of fewer than $k$ variables is contained in such a
subalgebra.  Hence the quotient satisfies every identity of its
specified variety on fewer than $k$ variables, yet fails $E_k$,
which holds in $A$ and all its subvarieties.  Choosing $k$ larger than any proposed bound on the number of variables
rules out a bounded-variable basis, and hence a finite basis.  The preceding two lemmas
give strictness, and local finiteness is inherited from $\mathcal U_2$.
The proper nonfinitely based subvarieties show that
$\mathcal L_2$ and $\mathcal L_3$ are not limit varieties.
\end{proof}

\section{Eventual recognition and strict separation}
\label{eventual:section}

We show that identities involving a fixed number of variables have
the same truth value in all sufficiently large finite recognition
quotients.  Their cofinite reduced product then generates an explicit
nonfinitely based variety with a decidable equational theory.  We also
show that it is properly contained in an explicitly axiomatized variety.

For $n\geq3$, retain $D_n^-$, $\Gamma$, $\Delta$, and the cap map
$c(v)_j=\min\{2,v_j\}$ from the punctured construction.  The dimension
of $\Gamma$ and $\Delta$ will always be specified by the surrounding
algebra.  Write
\[
 \sigma_n(S)\quad\Longleftrightarrow\quad S\subseteq\Delta
 \qquad(S\in D_n^-).
\]
The set $\Delta$ is closed under the modified closure $C_*$, so
\begin{equation}\label{eventual:conjunction}
 \sigma_n(S+T)\quad\Longleftrightarrow\quad
 \sigma_n(S)\ \text{and}\ \sigma_n(T).
\end{equation}
Define a relation on $D_n^-$ by
\begin{equation}\label{eventual:relation}
 S\sim_n T\quad\Longleftrightarrow\quad
 \begin{cases}
  \sigma_n(S)=\sigma_n(T),\\
  \sigma_n(SR)=\sigma_n(TR)&\text{for every }R\in D_n^-.
 \end{cases}
\end{equation}
Equality between predicates means agreement of their truth values.
The first line tests the elements without exterior multiplication;
no unit symbol is introduced.  Put
\[
 \mathcal Q_n=D_n^-/{\sim_n}.
\]
Thus $\mathcal Q_n$ denotes a quotient algebra, whereas $Q_n$ in the
punctured construction denotes a cube of vectors.

\begin{lemma}\label{eventual:quotients}
The relation $\sim_n$ is the largest congruence on $D_n^-$ that
preserves $\sigma_n$.  Its quotient $\mathcal Q_n$ is an explicitly
constructible finite algebra.  It satisfies the commutative
additively idempotent semiring identities, $x+xy=x$, {$x^2=x^3$, and $G_2$}.
It satisfies every identity of $\mathcal U_2$ involving fewer than
$n$ variables and fails the cyclic identity $E_n$ of
\eqref{nfb:cyclic}.
\end{lemma}

\begin{proof}
Lemma~\ref{descent:recognition}, applied to the $C_*$-closed set
$\Delta$, proves that $\sim_n$ is a congruence. If a congruence
$\rho$ preserves $\sigma_n$ and $S\mathrel\rho T$, then
$SR\mathrel\rho TR$ for every parameter $R$.  Applying preservation
of $\sigma_n$ to these pairs and to the original pair gives
$S\sim_n T$.  Hence $\rho\subseteq\sim_n$, as required.

By Lemma~\ref{nfb:quotient}, $D_n^-$ is a subalgebra of $D_n$,
so the stated semiring and power identities descend to the quotient.
{The recognition congruence used in the proof of
Lemma~\ref{nfb:quotient} is exactly $\sim_n$.  Thus
$\theta_n\subseteq\sim_n$, and every evaluation instance of
$G_2$ is identified in $\mathcal Q_n$.}
Lemma~\ref{nfb:smallvariables} gives all identities of $\mathcal U_2$
on fewer than $n$ variables. The principal assignments for $E_n$
give $\Gamma$ and $\Delta$, which the test without multiplication
separates. Hence $E_n$ fails in $\mathcal Q_n$.

The quotient can be constructed by enumerating the closed subsets
that define $D_n^-$, computing their operation tables, and identifying
elements whose finite lists of recognition values agree.  Closedness
is decidable by rational linear feasibility and the exceptional rule
for $\Delta$.  Thus finite enumeration and exact rational arithmetic
suffice.
\end{proof}

For a nonzero vector $p\in\{0,1,2\}^r$, let $x^p$ denote the
nonempty product of the factors $x_i^{p_i}$ with $p_i>0$.  Coordinates
with $p_i=0$ mean that the corresponding factor is absent; they do
not stand for constant terms.  For a nonempty set
$P\subseteq\{0,1,2\}^r\setminus\{0\}$, write
$t_P=\sum_{p\in P}x^p$.

\begin{lemma}[Principal-vector criterion]\label{eventual:principal}
Let $r\geq1$, let $P$ be a nonempty subset of
$\{0,1,2\}^r\setminus\{0\}$, and let
$q\in\{0,1,2\}^r\setminus\{0\}$.  The absorption
\begin{equation}\label{eventual:absorption}
 t_P+x^q=t_P
\end{equation}
fails in $\mathcal Q_n$ if and only if there are nonzero vectors
$a_1,\ldots,a_r\in\{0,1,2\}^n$ and a vector
$w\in\{0,1,2\}^n$ such that
\begin{equation}\label{eventual:principaltest}
 b_p=c\left(w+\sum_{i=1}^r p_i a_i\right)\in\Delta
 \quad(p\in P),\qquad
 b_q=c\left(w+\sum_{i=1}^r q_i a_i\right)\notin\Delta.
\end{equation}
The choice $w=0$ represents the test with no exterior multiplication.
It is notation in the vector calculation, not a constant term.
\end{lemma}

\begin{proof}
The semiring laws, $x+xy=x$, and $x^2=x^3$ give the Frobenius
identity.  Indeed, expansion and the absorption law give
\[
 (x+y)^3=x^3+x^2y+xy^2+y^3=x^2+y^2,
\]
while square stabilization identifies $(x+y)^3$ with $(x+y)^2$.
Induction extends this equality to the square of any finite nonempty
sum.  As a result, every capped monomial $x^p$ is additive separately
in every variable that occurs in it.  If a variable is absent, expanding over its possible choices only
repeats the same summand, which additive idempotence removes.

Every element $S\in D_n^-$ is the sum of the principal filters
$\mathord\uparrow a=C_*(\{a\})$ for its nonzero points $a\in S$.
Expand each assigned element in this fashion.  Expand an exterior
parameter $R\in D_n^-$ in the same fashion when that parameter is
present.  Diagonal expansion of squares shows that the raw vectors of
a monomial evaluation are precisely the vectors
$c(w+\sum_i p_i a_i)$, with independent choices of the indicated
principal points.  Closure is contained in $\Delta$ exactly when
every raw vector belongs to $\Delta$, because $\Delta$ is
$C_*$-closed.

Suppose \eqref{eventual:absorption} fails in $\mathcal Q_n$.
Some test in \eqref{eventual:relation} distinguishes its two
evaluations.  Distributivity and
\eqref{eventual:conjunction} show that the right evaluation has
predicate value true and that a raw vector contributed by $x^q$ on
the left has predicate value false.  Select principal points
producing that vector.  For variables absent from $x^q$, select any
point from their nonempty assigned elements.  Every raw evaluation of
every right monomial is in $\Delta$, so this same coupled selection
has $b_p\in\Delta$ for every $p\in P$.  If there is no exterior factor, set $w=0$; otherwise use the selected
nonzero principal point of that factor.  This proves necessity.

Conversely, assign $x_i$ to $\mathord\uparrow a_i$.  If $w\ne0$,
use the multiplicative test parameter $\mathord\uparrow w$; if
$w=0$, use the first line of \eqref{eventual:relation}.  By the principal product calculation and closure criterion above,
the right evaluation has predicate value true and the left evaluation
has predicate value false.  The two quotient classes are therefore
distinct, so the absorption fails.
\end{proof}

Recall from Lemma~\ref{nfb:recognition} that a vector belongs to
$\Delta$ precisely when it is positive with at least one coordinate
two, or it has exactly one zero coordinate whose cyclic predecessor
is two.  In the following two lemmas, a cyclic coordinate is deleted
or inserted simultaneously in every vector appearing in
\eqref{eventual:principaltest}.

\begin{lemma}[Padding]\label{eventual:padding}
Suppose \eqref{eventual:principaltest} is a failure witness in
dimension $n$, and some coordinate $j$ is positive in every right
output $b_p$, $p\in P$.  Then the absorption fails in
$\mathcal Q_N$ for every $N\geq n$.
\end{lemma}

\begin{proof}
Insert a coordinate immediately before $j$ in the cyclic order.
At this new coordinate put zero in every $a_i$ and put one in $w$.
The new coordinate of each output is then one.  The old coordinates
are unchanged.  Every assigned vector remains nonzero.  The new
exterior vector is nonzero, so it is a permitted test parameter even
when the original exterior context was absent.

A right output with no zero retains a coordinate two and thus remains
in $\Delta$.  For a right output with one zero, the insertion is
not immediately before that zero: the coordinate immediately after
the insertion is $j$, which is positive in every right output.
Its predecessor-two pair is consequently unchanged.  No new zero
has been introduced, so this output also remains in $\Delta$.

Consider the target $b_q$.  If it has at least two zeros, both
remain.  If it has exactly one zero, its predecessor is one because
the target is outside $\Delta$.  An insertion elsewhere leaves this
pair unchanged; an insertion immediately before the zero replaces
its predecessor by the new value one.  In either case the target
remains outside $\Delta$.  If it has no zero, it is the all-one
vector, and insertion of a further one preserves this property.
These are all possible target patterns.

The inserted coordinate is now positive in every right output, so
the construction can be repeated.  It produces a witness in every
larger dimension.
\end{proof}

\begin{lemma}[Compression]\label{eventual:compression}
Put $s=|P|$.  If \eqref{eventual:absorption} fails in
$\mathcal Q_n$ and $n>s$, there is a failure witness in dimension
$d$ with
\[
 3\leq d\leq r+2s+3
\]
having a coordinate positive in every right output.  In particular,
this witness can be padded to every dimension $N\geq d$.
\end{lemma}

\begin{proof}
Take the vectors supplied by Lemma~\ref{eventual:principal}.
Each right output has at most one zero.  The union of the zero
positions of all $s$ right outputs therefore has size at most $s$.
Choose a coordinate $j$ outside that union; it exists because $n>s$.

Retain a set of cyclic coordinates containing the following choices.
First retain one nonzero coordinate of each $a_i$, requiring at most
$r$ choices.  For each right output with a zero, retain that zero and
its immediate cyclic predecessor.  For each right output without a
zero, retain a coordinate where it is two.  This requires at most $2s$ coordinates.  If the target has at least two zeros, retain
two of them.  If it has exactly one zero, retain the zero and its
immediate predecessor.  If the target has no zero, it is all one and
needs no additional choice.  At most two coordinates are needed for the target.  Finally retain $j$, requiring at most one more.
Delete all other coordinates while preserving the cyclic order.
If the selected set has fewer than three coordinates, enlarge it
using original coordinates until it has three.  This is possible
because $n\geq3$, and the bound $r+2s+3\geq6$ is preserved.

Restriction commutes with the capped coordinatewise calculations.
Every restricted $a_i$ remains nonzero by the first choices.  A
right output with one zero retains that zero and its original
immediate predecessor, whose value is two.  These two coordinates
are still adjacent in the restricted cyclic order: no original
coordinate lay between them.  Every other retained coordinate of
this output is positive.  Thus it still has precisely one zero
with predecessor two.  A right output that originally had no zero
cannot acquire a zero by restriction and retains a chosen two.
Consequently no new zero or forbidden predecessor-one pair is
created in a right output.

If the target had at least two zeros, the retained two still exclude
it from $\Delta$.  If it had exactly one zero, the retained original
predecessor remains adjacent and has value one, preserving the
obstruction.  If it was all one, its restriction remains all one.
Adding original coordinates to the retained set does not affect these
conclusions: an originally adjacent predecessor pair remains adjacent,
and an originally positive coordinate cannot become zero.
The coordinate $j$ remains common and positive.  The restricted
exterior vector may become zero; in that case it denotes the
permitted test without an exterior multiplication.

The number $d$ of retained coordinates satisfies the required bound.
The absorption fails in $\mathcal Q_d$ by
Lemma~\ref{eventual:principal}, and it fails in every larger dimension
by Lemma~\ref{eventual:padding}.
\end{proof}

\begin{theorem}[Uniform finite-rank stabilization]
\label{eventual:stabilization}
For an absorption \eqref{eventual:absorption} in $r$ variables
with $|P|=s$, its truth value in $\mathcal Q_n$ is constant for
$n\geq r+2s+3$.  More generally, put
\[
 B(r)=r+2\cdot3^r+1\qquad(r\geq1).
\]
Every identity on at most $r$ variables has the same truth value in
all $\mathcal Q_n$ with $n\geq B(r)$.
\end{theorem}

\begin{proof}
Put $b=r+2s+3$.  A failure in dimension $b$ has a coordinate positive
in every right output, since $b>s$.  The padding lemma propagates it
to every dimension $n\geq b$.  Conversely, a failure in any
dimension $n\geq b$ compresses to a dimension $d\leq b$ with a
common positive coordinate.  Padding that witness gives a failure
in dimension $b$.  Thus failure in dimension $n$ is equivalent to failure in dimension
$b$ throughout this range.

Distributivity and exponent capping reduce every term on $r$
variables, in all the algebras $\mathcal Q_n$, to $t_P$ for a
nonempty subset
$P\subseteq\{0,1,2\}^r\setminus\{0\}$.  No empty sum or empty
product is used.  There are $3^r-1$ possible capped monomials.
Equality $t_P=t_T$ is equivalent to the conjunction of the
absorptions $t_P+x^q=t_P$ for $q\in T$ and
$t_T+x^p=t_T$ for $p\in P$.  Each absorption has at most
$3^r-1$ right monomials.  Since
\[
 r+2(3^r-1)+3=B(r),
\]
the first assertion supplies this single bound for every absorption
and hence every identity on at most $r$ variables.  The same bound therefore applies to any collection of such identities,
regardless of their lengths or the size of the collection.
\end{proof}

\begin{theorem}[An explicit reduced-product variety]
\label{eventual:main}
Define
\[
 \begin{gathered}
 M=\left(\prod_{n\geq3}\mathcal Q_n\right)/{\approx},
 \qquad \mathcal W_\infty=\mathsf V(M),\\
 a\approx b\quad\Longleftrightarrow\quad
 a_n=b_n\text{ for all sufficiently large }n.
 \end{gathered}
\]
Then $\mathcal W_\infty$ is a locally finite subvariety of
$\mathcal U_2$ with no identity basis using a bounded number of
variables.  For every identity $\eta$ on at most $r$ variables,
\begin{equation}\label{eventual:decision}
 \mathcal W_\infty\models\eta
 \quad\Longleftrightarrow\quad
 \mathcal Q_{B(r)}\models\eta.
\end{equation}
For every $n\geq B(r)$, $\mathcal Q_n$ satisfies all identities of
$\mathcal W_\infty$ on at most $r$ variables.  The equational theory
of $\mathcal W_\infty$ is decidable.
\end{theorem}

\begin{proof}
The relation $\approx$ is a congruence, since agreement of the inputs
outside finite sets of indices implies agreement of the outputs of
either coordinatewise operation outside the union of those sets.  An identity holds in $M$ if and only if it holds in
$\mathcal Q_n$ for all sufficiently large $n$.  For the forward implication, suppose, contrapositively, that the
identity fails at infinitely many indices.  Order each finite algebra
$\mathcal Q_n$ by the enumeration used in its construction.  At every failing
index choose the first witnessing tuple in the induced finite order;
at all other indices use the first available tuple.  The entries of
these tuples form elements of the direct product.  Their two term
evaluations differ at infinitely many indices and thus represent
distinct elements of $M$.  Conversely, if the identity holds at
all indices beyond some bound, every tuple of sequences has term
values agreeing beyond that bound, so the identity holds in $M$.

Let $\eta$ be any identity of $\mathcal U_2$, and let $r$ be the
number of its variables.  By Lemma~\ref{eventual:quotients},
$\mathcal Q_n\models\eta$ whenever $n>r$ and $n\geq3$.
The preceding equivalence gives $M\models\eta$.  Therefore
$M\in\mathcal U_2$ and $\mathcal W_\infty\subseteq\mathcal U_2$.
Every algebra satisfying the commutative additively idempotent
semiring laws and $x^2=x^3$ has at most
$2^{3^r-1}-1$ distinct term values generated by any fixed $r$
elements: distributivity and exponent capping give exactly the
nonempty polynomial representations used above as an upper bound.
This proves local finiteness.

The identities of an algebra and of its generated variety coincide.
Apply Theorem~\ref{eventual:stabilization} to the equivalence between
validity in $M$ and eventual validity in the factors.  This proves
\eqref{eventual:decision}.  Since $B(r)$ depends only on the number of variables, each
$\mathcal Q_n$, $n\geq B(r)$, satisfies every $r$-variable identity
of $\mathcal W_\infty$.  To decide
an identity, count its variables, construct the finite algebra
$\mathcal Q_{B(r)}$ by Lemma~\ref{eventual:quotients}, and test
all its assignments.  This procedure terminates, although no bound on its efficiency is
claimed.

Suppose an identity basis $\Sigma$ of $\mathcal W_\infty$ used
at most $r$ variables in every member, after renaming variables
where necessary.  Take $n\geq B(r)$.  The simultaneous satisfaction
statement gives $\mathcal Q_n\models\Sigma$.  By
Lemma~\ref{eventual:quotients}, however, $\mathcal Q_n$ fails $E_n$.
This cyclic identity holds in $A$, hence in $\mathcal U_2$ and in
$\mathcal W_\infty$.  Soundness of equational deduction contradicts
the assumption that $\Sigma$ is a basis.  Thus no bounded-variable
basis exists; in particular no finite basis exists.
\end{proof}

\subsection{Separation from the safe-identity variety}
\label{separation:section}

Let $\Sigma_{\mathrm{safe}}$ be the explicitly specified set of all
absorptions satisfying the hypotheses of Lemma~\ref{descent:safe}:
the target is squarefree, every right-hand exponent is at most two,
the intersection of the doubled-variable sets is empty, and some
right-hand monomial omits at least two target variables.  All
monomials are nonempty and use only target variables.  Define
\[
 \mathcal W_*=
 \mathcal U_2\cap\operatorname{Mod}
             (\{G_2\}\cup\Sigma_{\mathrm{safe}}).
\]
This gives an explicit infinite equational presentation without a
maximality choice.  We next show that it defines a variety properly
containing the eventual recognition variety.

\begin{theorem}\label{separation:strict}
For the identity
\begin{equation}\label{separation:eta4}
 \eta_4:\quad
 xyzw+xy^2z^2+xy^2w^2+xz^2w^2
       =xy^2z^2+xy^2w^2+xz^2w^2,
\end{equation}
every recognition quotient $\mathcal Q_k$, $k\geq3$, satisfies
$\eta_4$, but some finite member of $\mathcal W_*$ does not.
Consequently
\[
 \mathcal W_\infty\subsetneq\mathcal W_*\subseteq\mathcal L_4.
\]
In particular $\mathcal L_4$ is not a limit variety.
\end{theorem}

\begin{proof}
By Lemma~\ref{descent:safe}, every safe absorption holds in every
$\mathcal Q_k$.  The quotients also satisfy $G_2$ by
Lemma~\ref{nfb:quotient}.  The reduced-product characterization in
Theorem~\ref{eventual:main} therefore gives
$\mathcal W_\infty\subseteq\mathcal W_*$.  Both $\alpha$ and
$\beta$ belong to $\Sigma_{\mathrm{safe}}$, so
$\mathcal W_*\subseteq\mathcal L_4$.

We prove $\eta_4$ in $\mathcal Q_k$ using the principal-vector
criterion of Lemma~\ref{eventual:principal}.  Take an assignment by
four nonzero principal vectors and an exterior multiplicative
context.  Let $a$ be the capped target vector for $xyzw$.  If $a$
has at least two zeros, all right-hand vectors retain both zeros and
lie outside $\Delta$.  If $a$ has one zero at $j$ and predecessor
value one, this value comes from the context alone or from one
unique assigned variable.  The doubled-variable sets on the right
are $\{y,z\},\{y,w\},\{z,w\}$, with empty intersection.
Thus at least one right-hand vector retains the zero at $j$ and has
predecessor value at most one, placing it outside $\Delta$.

The remaining bad case is $a=\delta$.  The four assigned nonzero
vectors then occupy disjoint nonempty supports
$B_x,B_y,B_z,B_w$, and the context has disjoint support.  The three
right monomials respectively omit $w,z,y$.  For all their vectors
to belong to $\Delta$, the blocks $B_y,B_z,B_w$ must be singletons.
The predecessor of each corresponding zero coordinate must have
value two in its right monomial.  It cannot belong to $B_x$ or to
the context, since those coordinates have value one in every
right-hand vector.  Hence the three-point set
$B_y\cup B_z\cup B_w$ is closed under the immediate cyclic
predecessor map.  A nonempty subset closed under this map is the
whole coordinate cycle: repeated predecessors of one point traverse
every position.  This contradicts the additional nonempty block
$B_x$.  When $k=3$, already four disjoint nonempty blocks are
impossible.  The required recognition implication holds for all
$k\geq3$, and Lemma~\ref{descent:recognition} proves $\eta_4$ in
every $\mathcal Q_k$, hence in $\mathcal W_\infty$.

For strictness, work in the ordinary capped convex free algebra
$C_{2,4}$ of Proposition~\ref{cap:free}.  Fix coordinate sets
$U=\{x,y,z,w\}$ and $T_0=\{y,z,w\}$, and put
$\delta=(1,1,1,1)$.  In its integer cube define
\begin{align*}
 K={}&\bigl(\mathord{\uparrow}\delta\setminus\{\delta\}\bigr)
       \cup\mathord{\uparrow}(1,0,2,2)\\
    &\cup\mathord{\uparrow}(1,2,0,2)
       \cup\mathord{\uparrow}(1,2,2,0).
\end{align*}
Thus $K$ consists of the positive nonlinear vectors and the vectors
having exactly one zero at $i\in T_0$, with $x$ positive and both
other $T_0$-coordinates equal to two.

The set $K$ is convex upward closed in the integer cube.  Every point
has weight at least five, so no convex combination lies below
$\delta$, of weight four.  If an integer vector above a convex
combination has a zero at $i\in T_0$, every point with positive
coefficient belongs to the unique displayed principal upset with
zero there.  Their average and the vector above it remain in that
upset.  A zero at $x$ is impossible, since $x$ is positive throughout
$K$.  Every remaining positive integer vector belongs to $K$.

For a nonlinear vector of fixed support $T$ and residual set $J$,
membership in $K$ means either $T=U$, or
$T=U\setminus\{i\}$ for $i\in T_0$ and $J\subseteq\{x\}$.
All other supports are excluded.  In a principal $G_2$ instance,
support is fixed and the two right residual sets have union equal
to the left residual set, as in \eqref{descent:residualunion}.
This criterion proves $G_2$ invariance for recognition by $K$ in all
contexts.

We check every member of $\Sigma_{\mathrm{safe}}$, allowing any
number of its variables.  Let $a$ be a bad target point in a common
principal assignment and context.  If $a$ has at least two zeros,
all right points retain them.  If it has a zero at $x$, all right
points retain that zero.  The remaining bad point with a zero has
exactly one zero at $i\in T_0$ and value one at another
$T_0$-coordinate $j$.  That one comes from the context or from one
unique variable.  In the context case every right point retains it;
otherwise the empty intersection of doubled-variable sets supplies
a right point with value at most one at $j$.  It retains zero $i$
and is outside $K$.  The last bad possibility is $a=\delta$.
A right monomial omitting at least two variables omits two disjoint
nonempty assignment blocks, so has at least two zeros and lies
outside $K$.  This covers every point outside $K$.  By
Lemma~\ref{descent:recognition}, every instance of a safe absorption,
and every instance of $G_2$, belongs to $\equiv_K$.

The finite recognition quotient $C_{2,4}/{\equiv_K}$ belongs to
$\mathcal U_2$ and satisfies $G_2$ and $\Sigma_{\mathrm{safe}}$.
It therefore belongs to $\mathcal W_*$.  At the principal generators
the right-side exponent vectors of $\eta_4$ lie in $K$, while
$\delta$ does not.  The identity context separates the two quotient
values, so $\eta_4$ fails there.  This proves strict inclusion.

Finally, Theorem~\ref{eventual:main} makes
$\mathcal W_\infty$ nonfinitely based.  It is proper in
$\mathcal L_4$, since its separating identity fails in a member of
$\mathcal W_*\subseteq\mathcal L_4$.  Hence $\mathcal L_4$ has a proper nonfinitely based subvariety, so
it is not a limit variety.
\end{proof}
\section{Cycle quotients and subvariety intervals}\label{card:section}
We now require both neighbors of a zero coordinate to have value two.
The resulting finite quotients yield a continuum-sized subvariety
interval through an independent family of identities. Their reduced
product generates a proper nonfinitely based subvariety of
$\mathcal W_\infty$.

\subsection{Recognition quotients attached to cycles}
For $n\geq3$, use the capped convex algebra
$C_n=C_{2,n}$ from Section~\ref{sec:convex}, on
$T_n=\{0,1,2\}^n\setminus\{\boldsymbol0\}$, and write
$c=c_2$ and $\operatorname{cl}$ for its cap map and closure.

\begin{lemma}\label{card:cappedmodel}
The algebra $C_n$ is a finite member of $\mathcal U_2$.
Each element is a finite sum of principal upper sets, and
\[
 S^d=\operatorname{cl}\{c(ds):s\in S\}\quad(d\geq1),
 \qquad (S+T)^2=S^2+T^2.
\]
\end{lemma}
\begin{proof}
These are Proposition~\ref{cap:membership} and
Remark~\ref{rem:principal} with $m=2$.
\end{proof}

Read indices modulo $n$, and let $\delta_n=(1,\ldots,1)$.
For each $i$, define $d_{n,i}$ to have coordinate zero at $i$,
coordinate two at the two cyclic neighbors $i-1,i+1$, and
coordinate one elsewhere.  Its neighbors are distinct because
$n\geq3$.  For upward closures relative to $T_n$, put
\[
 K_n=\bigl(\mathord{\uparrow}\delta_n\setminus\{\delta_n\}\bigr)
       \cup\bigcup_{i=1}^n\mathord{\uparrow}d_{n,i}.
\]

\begin{lemma}\label{card:recognitionset}
The set $K_n$ is nonempty and $\operatorname{cl}$-closed.  A vector belongs to it
exactly when either it is positive and differs from $\delta_n$,
or it has one zero coordinate whose two cyclic neighbors have value
two and whose other coordinates are positive.
\end{lemma}

\begin{proof}
The description of membership follows from the definition as a union
of upsets. Every $d_{n,i}$ has weight $n+1$, and every positive
vector different from $\delta_n$ has weight at least $n+1$.
Let $a$ be an integer vector above a convex combination of points
of $K_n$.
If $a$ is positive, its only possible exclusion is $\delta_n$,
whose weight $n$ is too small.  If $a_i=0$, every point with
positive coefficient in the combination must also have zero at
$i$.  Such points belong to the single principal upset
$\mathord{\uparrow}d_{n,i}$; all other pieces of $K_n$ are
positive there.  Their average and the vector above it are at
least $d_{n,i}$, placing $a$ in $K_n$.  This also excludes a
second zero.  The required closedness follows.
\end{proof}

For $S\in C_n$, write $\sigma_n(S)$ for $S\subseteq K_n$.
Define $S\sim_n T$ when
\[
 \sigma_n(SR)=\sigma_n(TR)\quad\text{for every }R\in C_n,
 \qquad \sigma_n(S)=\sigma_n(T),
\]
where equality denotes agreement of truth values.  The last test
is the context without multiplication; it does not adjoin a
multiplicative constant.  Set
\[
 S_n=C_n/{\sim_n}.
\]

\begin{lemma}\label{card:quotient}
The relation $\sim_n$ is a semiring congruence.  Thus $S_n$ is a
finite member of $\mathcal U_2$, and in particular of $\mathcal V$.
\end{lemma}

\begin{proof}
Apply Lemma~\ref{descent:recognition} to the closed set $K_n$
of Lemma~\ref{card:recognitionset}. This gives the congruence
$\sim_n$. Since $C_n\in\mathcal U_2$ by
Lemma~\ref{card:cappedmodel}, the finite quotient belongs to
$\mathcal U_2$ as well.
\end{proof}

For each $n$, the quotient can be constructed by enumerating the
closed subsets of $T_n$, computing the operation tables, and recording
$\sigma_n(S)$ and $\sigma_n(SR)$ for all parameters.
Two elements are identified precisely when their rows agree.
Thus the construction requires no maximality choice.

\begin{lemma}[Principal-vector criterion]\label{card:principal}
A capped absorption $P+q=P$ on $r$ variables fails in $S_n$ if
and only if there are nonzero vectors $a_1,\ldots,a_r\in T_n$
and $w\in\{0,1,2\}^n$ such that
\[
 b_p=c\Bigl(w+\sum_jp_ja_j\Bigr)\in K_n\quad(p\in P),
 \qquad
 b_q=c\Bigl(w+\sum_jq_ja_j\Bigr)\notin K_n.
\]
The parameter value $w=0$ means omission of the exterior factor.
\end{lemma}

\begin{proof}
Write the assigned closed sets and the exterior parameter as sums
of principal upper sets. By Lemma~\ref{card:cappedmodel}, squares
expand diagonally. Since $K_n$ is closed, a polynomial value passes
a recognition test if and only if all its raw generating vectors
belong to $K_n$. Thus failure of the absorption gives a target vector outside $K_n$
with all right vectors inside it. Choose the principal points producing
that target, and arbitrary points for variables absent from it.
Using these choices in every right monomial gives the displayed
conditions. Conversely, these conditions separate the two values
at the principal assignments. The zero exterior vector represents
the test without a factor.
\end{proof}

\subsection{An independent cycle family of identities}

For $r\geq3$ put
\[
 q_r=\prod_{j=1}^r x_j,\qquad
 p_{r,i}=\prod_{j\ne i}x_j^{\,1+\mathbf1_{\{i-1,i+1\}}(j)},
\]
where the indices are read modulo $r$. The product defining
$p_{r,i}$ {has exactly $r-1\geq2$ variables}, each with exponent
one or two. Thus no empty product or zeroth power occurs. Define
\begin{equation}\label{card:eta}
 \eta_r:\qquad q_r+\sum_{i=1}^r p_{r,i}
                      =\sum_{i=1}^r p_{r,i}.
\end{equation}
Its right-side exponent vectors are the $d_{r,i}$.

\begin{proposition}\label{card:independence}
For every $n,r\geq3$,
\[
 S_n\models\eta_r\quad\Longleftrightarrow\quad n\ne r.
\]
\end{proposition}

\begin{proof}
First suppose $n=r$, and assign $x_i$ to the class of the principal
generator $\operatorname{cl}(\{e_i\})$. All right-side generating
vectors $d_{n,i}$ belong to $K_n$, so their closure is contained
there. The left value contains $\delta_n\notin K_n$.
The test without an exterior factor therefore distinguishes the two
values, and $\eta_n$ fails in $S_n$.

Now suppose $n\ne r$ and, towards a contradiction, take a principal
witness as in Lemma~\ref{card:principal}. Write $w$ for the
context vector and put
\[
 a=c(w+a_1+\cdots+a_r),\qquad
 b_i=c\Bigl(w+\sum_{j\ne i}
          \bigl(1+\mathbf1_{\{i-1,i+1\}}(j)\bigr)a_j\Bigr).
\]
All $b_i$ belong to $K_n$, whereas $a$ does not. We consider
the possible forms of the target given by
Lemma~\ref{card:recognitionset}.

If $a$ has at least two zeros, all assigned vectors and $w$ vanish
there.  Every $b_i$ retains both zeros and is outside $K_n$,
a contradiction.

Suppose $a$ has exactly one zero at position $s$.  Since $a$ is
outside $K_n$, one cyclic neighbor $t$ of $s$ has value one.
The uncapped sum at $t$ is one.  If it comes from $w$, each $b_i$
retains value one at $t$ and zero at $s$, and is outside $K_n$.
Otherwise it comes from a unique vector $a_j$, with all other
vectors and $w$ zero there.  The monomial $p_{r,j}$ omits that
variable, so $b_j$ has zeros at both $s,t$ and is outside $K_n$.
This is again a contradiction.

The remaining bad target is $a=\delta_n$.  No capping occurs in
a sum with this value.  The assigned nonzero vectors are zero-one
vectors on pairwise disjoint nonempty supports
$B_1,\ldots,B_r$, while $w$ is supported on the disjoint remaining
set $B_0$.  The vector $b_i$ has zeros exactly on $B_i$.
Membership in $K_n$ forces every $B_i$ to be a singleton.  Write
$B_i=\{f(i)\}$, obtaining an injection from the $r$-cycle into
the coordinate set of the $n$-cycle.

In $b_i$, the two coordinates with value two are precisely
$f(i-1),f(i+1)$.  Its zero is at $f(i)$.
Lemma~\ref{card:recognitionset} requires both neighbors of $f(i)$
in the $n$-cycle to be in this two-point set.  Both sets have
cardinality two, so
\[
 \{f(i)-1,f(i)+1\}=\{f(i-1),f(i+1)\}.
\]
The nonempty image of $f$ is therefore closed under cyclic
neighbors. Taking successive neighbors shows that the image contains
the entire $n$-cycle. Thus the injection is surjective and $r=n$,
a contradiction. Hence no principal witness exists when $r\ne n$.
\end{proof}

\subsection{The cardinality theorem}

\begin{theorem}\label{card:main}
The variety $\V(\N,\min,+_{\rm ar})$ has exactly
$2^{\aleph_0}$ subvarieties.  Already its proper locally finite
subvariety $\mathcal U_2$ has exactly $2^{\aleph_0}$ subvarieties,
and $2^{\aleph_0}$ of these are nonfinitely based in the fixed
constant-free language.
\end{theorem}

\begin{proof}
For each $I\subseteq\{4,5,6,\ldots\}$ define
\[
 \mathcal V_I=\V\bigl(\{S_3\}\cup\{S_n:n\in I\}\bigr).
\]
Every generator belongs to $\mathcal U_2$, so
$\mathcal V_I\subseteq\mathcal U_2\subseteq\mathcal V$.
For $m\geq4$, Proposition~\ref{card:independence} gives
\[
 \mathcal V_I\models\eta_m
       \quad\Longleftrightarrow\quad m\notin I.
\]
Indeed, if $m\in I$, the identity fails in $S_m$. If $m\notin I$,
it holds in all the generators and hence in the variety they generate.
Distinct subsets $I$ therefore give distinct subvarieties, and there
are $2^{\aleph_0}$ such subsets.

For the upper bound, fix a countable set of variables. There are
countably many terms and identities in a finite signature. Since
every variety is determined by its identities, there are at most
$2^{\aleph_0}$ subvarieties. This proves both cardinality assertions.

The identity $x^2=x^3$ fails in $A$ at the natural number one, so
$\mathcal U_2$ is proper.  Every term on $r$ generators reduces
to a nonempty polynomial in the $3^r-1$ nonzero capped monomials.
There are only finitely many such polynomials; thus every finitely
generated free algebra of $\mathcal U_2$ is finite.  Every finitely
generated algebra is a homomorphic image of the corresponding free
algebra, proving local finiteness.  Finally the countable identity syntax has only
countably many finite subsets, so there are only countably many
finitely based varieties.  Removing that countable set from the
continuum family leaves $2^{\aleph_0}$ nonfinitely based
subvarieties, all locally finite.
\end{proof}

\subsection{The cycle reduced product}
\label{smaller:section}

We continue with the finite ordinary recognition quotients $S_n$,
the closed sets $K_n$, and the identities $\eta_n$ used above.
We have $S_n\in\mathcal U_2$ and
\begin{equation}\label{smaller:independence}
 S_n\models\eta_m\quad\Longleftrightarrow\quad n\ne m
 \qquad(n,m\geq3)
\end{equation}
by Proposition~\ref{card:independence}.  Write $\mathcal Q_n$ for
the punctured recognition quotients of
Lemma~\ref{eventual:quotients}.  Recall the difference between the two constructions:
$S_n$ requires \emph{both} neighbors of a zero to have value two,
whereas $\mathcal Q_n$ requires only its predecessor to have value two.

\begin{lemma}\label{smaller:comparison}
For every $n\geq3$,
\[
 \operatorname{Id}(\mathcal Q_n)
 \subseteq\operatorname{Id}(S_n).
\]
Consequently $S_n$ satisfies every identity of
$\mathcal W_\infty$ on at most $r$ variables whenever
$n\geq B(r)=r+2\cdot3^r+1$.
\end{lemma}

\begin{proof}
Since every identity reduces to capped monomial absorptions, it
suffices to transfer a failed absorption from $S_n$ to
$\mathcal Q_n$. By Lemma~\ref{card:principal}, choose nonzero
principal vectors and an exterior vector such that every right output
belongs to $K_n$ but the target does not. As before, a zero exterior
vector represents the test without an extra factor.
The inclusion $K_n\subseteq\Delta_n$ puts every right output
in $\Delta_n$.

A target with two zeros is outside $\Delta_n$.  A target with no
zero that is outside $K_n$ is $\delta_n$, also outside
$\Delta_n$.  The remaining target has exactly one zero and a
neighbor one.  If that neighbor is the predecessor, it is outside
$\Delta_n$.  Otherwise reflect the entire cyclic coordinate order.
Reflection interchanges predecessor and successor and commutes with
the coordinatewise capped calculations. The assigned vectors remain
nonzero, and $K_n$ is preserved because its definition treats both
neighbors equally. The reflected right outputs therefore still belong
to $K_n$, while the reflected target has predecessor one and lies
outside $\Delta_n$.  Lemma~\ref{eventual:principal} now supplies a failure
in $\mathcal Q_n$.  This proves the identity inclusion.  The last
assertion follows from Theorem~\ref{eventual:main}.
\end{proof}

\begin{lemma}[Padding and compression]\label{smaller:compression}
Consider an absorption $t_P+x^q=t_P$ on $r\geq1$ variables, with
$s=|P|\geq1$ capped right monomials.
\begin{enumerate}
\item If a principal failure witness in $S_n$ has two consecutive
coordinates positive in every right output, it extends to a failure
witness in $S_N$ for every $N\geq n$.
\item If the absorption fails in $S_n$ and $n>2s$, it has a failure
witness in some dimension $d$ with
\[
 3\leq d\leq r+3s+4,
\]
having two consecutive coordinates positive in every right output.
\end{enumerate}
\end{lemma}

\begin{proof}
For the first assertion, let $u,v$ be consecutive coordinates with
the stated property. Insert a coordinate between them, assigning zero
to every variable vector and one to the exterior vector there.
Every output gains a coordinate of value one. The assigned vectors
remain nonzero, and the new exterior vector is nonzero even if
the original test had no exterior factor.

A positive right output retains a coordinate two.  A right output
with one zero has neither $u$ nor $v$ as that zero, so insertion
does not change either of its two neighboring coordinates.  Both
remain two.  No new zero appears, and every right output remains
in the larger recognized set.  A bad target with two zeros retains
them, and an all-one target remains all one.  A target with one
zero and a neighbor one remains bad as well: insertion elsewhere
preserves its neighbors, while insertion adjacent to its zero gives
it a new neighbor one.  The new coordinate and $u$ are consecutive and positive in every
right output, so the insertion can be repeated. This proves the
first assertion.

For the second, let $Z$ be the union of the zero positions of the
right outputs.  Each such output has at most one zero, so $|Z|\leq s$.
There are consecutive positions $u,v$ outside $Z$.  Otherwise the
successor of every position outside $Z$ would lie in $Z$, giving
an injection from the complement of $Z$ into $Z$ and forcing
$n\leq2|Z|\leq2s$.

Retain one nonzero coordinate of each assigned vector, requiring at
most $r$ coordinates.  For each right output with a zero, retain
that zero and its two original immediate neighbors.  For each
positive right output, retain a coordinate two.  These choices use
at most $3s$ coordinates.  For a target with at least two zeros,
retain two of them; for a target with one zero, retain that zero
and an original neighboring coordinate one; an all-one target needs
no extra choice.  This uses at most two more coordinates.  Retain
also $u,v$, requiring at most two.  Delete the other coordinates
while maintaining cyclic order.  If necessary retain further original
coordinates until at least three remain; this is possible since
$n\geq3$, and the upper bound $r+3s+4\geq8$ is preserved.

Restriction commutes with the capped vector calculations and leaves
the assigned vectors nonzero. Every positive right output retains
a coordinate two. Each retained zero of a right output keeps both
of its original neighbors, still adjacent to it. These neighbors
remain distinct because at least three coordinates have been retained.  Every other coordinate of this
output is positive, so it has precisely one zero and both neighbors
two.  The target retains either two zeros, a zero with an adjacent
one, or the all-one pattern.  The original consecutive pair $u,v$
remains consecutive and common positive.  If further original coordinates are retained to reach dimension three,
these conclusions still hold: each original adjacent pair stays adjacent,
and all other coordinates of a one-zero output were positive.
If restriction makes the exterior vector zero, omit the exterior factor.
This gives a witness of the required dimension with the stated properties.
\end{proof}

\begin{theorem}[Uniform stabilization of the ordinary quotients]
\label{smaller:stabilization}
For an absorption on $r$ variables with $s$ capped right monomials,
its truth value in $S_n$ is constant for $n\geq r+3s+4$.
More generally, every identity on at most $r$ variables has the same
truth value in all $S_n$ with
\[
 D(r)=r+3\cdot3^r+1\leq n.
\]
\end{theorem}

\begin{proof}
Put $b=r+3s+4$, which is greater than $2s$.  A failure at $b$
has a consecutive pair positive in all right outputs, by the counting
argument in Lemma~\ref{smaller:compression}.  Padding propagates
it to every $n\geq b$.  Conversely, a failure at any $n\geq b$
compresses to some $d\leq b$ with such a pair, and padding returns
a failure at $b$.  These implications prove the first assertion.

There are $3^r-1$ nonzero capped monomials on $r$ variables.
Under the common semiring laws and square stabilization, every identity
is equivalent to the monomial absorptions between its two capped
polynomials in both directions. Each absorption has at most $3^r-1$
right monomials. Hence the bound
\[
 r+3(3^r-1)+4=D(r)
\]
applies to all identities on at most $r$ variables, irrespective
of their lengths or number.
\end{proof}

To distinguish the two eventual varieties, we use the following
identities. For $k\geq3$, let
\[
 d_k=t_1\cdots t_k,\qquad
 g_{k,i}=t_i^2\prod_{j\notin\{i,i+1\}}t_j,\qquad
 \rho_k=\sum_{i=1}^k d_k t_i,
\]
with cyclic indices.  All monomials are nonempty.  Write
\begin{equation}\label{smaller:pathidentity}
 \mathsf C_k:\quad
 d_k+\sum_{i=1}^{k-1}g_{k,i}+\rho_k
   =\sum_{i=1}^{k-1}g_{k,i}+\rho_k.
\end{equation}
The symbol $\mathsf C_k$ names an identity, not a convex algebra.

\begin{lemma}\label{smaller:path}
Every $S_n$, $n\geq3$, satisfies every $\mathsf C_k$, $k\geq3$,
and satisfies $G_2:z^2xy=z^2x^2y+z^2xy^2$.
The identity $\mathsf C_3$ fails in $\mathcal W_\infty$.
\end{lemma}

\begin{proof}
First consider $\mathsf C_k$. Take a principal evaluation with
exterior vector $w$ and suppose that its target
$a=c_2(w+a_1+\cdots+a_k)$ lies outside $K_n$.
If it has at least two zeros, every right output retains them.
If it has one zero $h$ and a neighbor $l$ of value one, the one
comes from the context or a unique assigned vector $a_j$.  In the
first case every right output retains value one at $l$.  In the
second case choose a saturation term $d_k t_i$ with $i\ne j$;
its output retains the zero at $h$ and value one at $l$.  Thus
some right output is outside $K_n$.

If the target is $\delta_n$, the assigned nonzero variables have
pairwise disjoint nonempty supports.  The monomials
$g_{k,1},\ldots,g_{k,k-1}$ omit respectively the supports of
$a_2,\ldots,a_k$.  Membership of all their outputs in $K_n$
would force each of those supports to be a singleton.  The term
$g_{k,2}$ is present because $k\geq3$.  It doubles only $a_2$
and omits $a_3$.  {Since $a_2$ has singleton support,
its output has exactly one coordinate}
two and a zero at the position of $a_3$.  It cannot have two distinct neighbors of value two and hence lies
outside $K_n$. Thus, in every principal context, a failed target
forces a failed right output. By Lemma~\ref{card:principal},
$\mathsf C_k$ holds in $S_n$.

For $G_2$, use principal nonzero vectors $x,y,z$ and an exterior
vector $w$, and put
\[
 a=c_2(x+y+2z+w),\quad
 b=c_2(a+x),\quad
 d=c_2(a+y).
\]
Upward closure gives $b,d\in K_n$ if $a\in K_n$.  A bad target
$a$ cannot be $\delta_n$, because nonzero $z$ supplies a two.
Two zeros persist in both $b,d$.  In the remaining bad case $a$
has one zero and a neighbor whose uncapped sum is one.  At most one
of $x,y$ is positive there, so at least one of $b,d$ retains that
neighbor one while retaining the zero.  Thus $a\in K_n$ if and
only if $b,d\in K_n$.  Principal expansion, including all exterior
tests, proves $G_2$ in the quotient.

Finally, in each $\mathcal Q_n$, $n\geq3$, assign $t_1,t_2,t_3$
the principal vectors $e_1,e_2,e_3$.  Put the exterior vector equal
to one at the remaining coordinates and zero at the first three;
when $n=3$ use no exterior factor.  The target is $\delta_n$.
The output of $g_{3,1}$ has its only zero at coordinate two and
predecessor value two at coordinate one.  The output of $g_{3,2}$
has its only zero at coordinate three and predecessor value two at
coordinate two.  Every saturation output is positive with a two.
Thus every right output belongs to $\Delta_n$ and the target does
not.  The principal criterion refutes $\mathsf C_3$ in every
$\mathcal Q_n$.  The reduced-product criterion in
Theorem~\ref{eventual:main} proves its failure in
$\mathcal W_\infty$.
\end{proof}

\begin{theorem}\label{smaller:main}
Let
\[
 \begin{gathered}
 N=\left(\prod_{n\geq3}S_n\right)/{\approx},
 \qquad\mathcal Y=\mathsf V(N),\\
 a\approx b\quad\Longleftrightarrow\quad
 a_n=b_n\text{ for all sufficiently large }n.
 \end{gathered}
\]
Then $\mathcal Y$ is a proper locally finite subvariety of
$\mathcal W_\infty$ with no identity basis using a bounded number
of variables.  Its equational theory is decidable, and for every
identity $\theta$ on at most $r$ variables,
\begin{equation}\label{smaller:decision}
 \mathcal Y\models\theta
 \quad\Longleftrightarrow\quad S_{D(r)}\models\theta.
\end{equation}
Every $S_n$ with $n\geq D(r)$ satisfies all identities of
$\mathcal Y$ on at most $r$ variables.  In particular,
$\mathcal W_\infty$ is not a limit variety.
\end{theorem}

\begin{proof}
The reduced-product argument in the proof of
Theorem~\ref{eventual:main} applies to $(S_n)$: an identity holds
in $N$ exactly when it holds in all sufficiently large factors.

Every identity of $\mathcal W_\infty$ uses finitely many variables.
By Lemma~\ref{smaller:comparison}, it holds in $S_n$ for all
sufficiently large $n$, and hence in $N$. Therefore
$\mathcal Y\subseteq\mathcal W_\infty$. Containment in
$\mathcal U_2$ also gives local finiteness.  Combining eventual validity
with Theorem~\ref{smaller:stabilization} proves
\eqref{smaller:decision} and its simultaneous assertion.  Exact
construction of the finite algebra $S_{D(r)}$ followed by testing
all assignments gives a terminating decision procedure.

By \eqref{smaller:independence}, each $\eta_m$ fails only in
the factor $S_m$, so it holds in $N$ and $\mathcal Y$.
Suppose that $\mathcal Y$ had a basis using at most $r$ variables
in each identity. For $n\geq D(r)$, the simultaneous assertion
would make $S_n$ satisfy the whole basis. But it fails $\eta_n$,
an identity of $\mathcal Y$, contradicting soundness of equational
consequence. Hence there is no bounded-variable basis, and in
particular no finite basis.

By Lemma~\ref{smaller:path}, $\mathsf C_3$ holds in every $S_n$,
and hence in $N$ and $\mathcal Y$, but fails in
$\mathcal W_\infty$. Thus the inclusion is proper, and this
proper nonfinitely based subvariety shows that $\mathcal W_\infty$
is not a limit variety.
\end{proof}

\begin{proposition}\label{smaller:upperinterval}
The interval $[\mathcal W_\infty,\mathcal U_2]$ contains continuum
many distinct nonfinitely based subvarieties.  More generally, every
variety in $[\mathcal W_\infty,\mathcal V]$ has no identity basis
using a bounded number of variables.
\end{proposition}

\begin{proof}
First every $\eta_m$ holds in $\mathcal W_\infty$.  To prove
this, test it in $\mathcal Q_n$ for $n>m$ using principal vectors.
If its target has two zeros, every right output retains them.  If
the target has one zero and predecessor one, that one comes from
the context or a unique assigned variable.  In the first case every
right output retains predecessor one; in the second, the monomial
omitting that variable has predecessor zero.  Thus a right output
is outside $\Delta_n$.
For an all-one target, the assigned nonzero variables occupy disjoint
nonempty blocks.  Each is omitted by a right monomial, so membership
of all right outputs in $\Delta_n$ forces the blocks to be
singletons, at distinct positions $h_i$.  The physical predecessor
of $h_i$ must occur at one of the doubled variables of that right
monomial, hence in the set of positions $\{h_i\}$.  This nonempty
set is predecessor-closed in the physical coordinate cycle and
therefore consists of every coordinate.  It would give $m=n$,
contrary to $n>m$.  These cases establish the principal implication
and hence the identity in every $\mathcal Q_n$ for $n>m$.
The reduced-product criterion proves it in $\mathcal W_\infty$.

For $X\subseteq\{3,4,\ldots\}$ put
\[
 \mathcal Z_X=\mathcal W_\infty\vee
 \mathsf V(\{S_n:n\in X\}).
\]
For the empty set the second joinand means the trivial variety.
The preceding paragraph and \eqref{smaller:independence} give
\[
 \mathcal Z_X\models\eta_m\quad\Longleftrightarrow\quad m\notin X.
\]
Thus these varieties are pairwise distinct and give continuum many
members of the stated interval.

Let $\mathcal W_\infty\subseteq\mathcal W\subseteq\mathcal V$,
and suppose a basis of $\mathcal W$ used at most $r$ variables.
Every identity in that basis holds in $\mathcal W_\infty$.
For $n\geq B(r)$, Theorem~\ref{eventual:main} puts every such
identity in $\mathcal Q_n$.  But $\mathcal Q_n$ fails $E_n$,
which holds in $\mathcal V$ and hence in $\mathcal W$.
This contradicts soundness of equational consequence. The last
assertion follows, and in particular every $\mathcal Z_X$ is
nonfinitely based.
\end{proof}

\section{Symmetric recognition quotients}
\label{clique:section}

We now replace the condition on the two neighbors of a zero by a
condition on every nonzero coordinate of the vector. For $n\geq3$, put
\[
 H_n=\{v\in\{0,1,2\}^n:v\text{ has all coordinates positive},\ v\ne\delta_n\}
       \ \cup\ \{p'_{n,i}:1\leq i\leq n\},
\]
where $\delta_n=(1,\ldots,1)$ and $p'_{n,i}$ has zero at $i$ and
two at every other coordinate.  Equivalently, $H_n$ contains the
positive non-all-one vectors and the vectors with exactly one zero
and all other coordinates two.  Each $p'_{n,i}$ is maximal among
cube vectors with zero at $i$; increasing its zero puts it in the
positive part of $H_n$.

\begin{lemma}\label{clique:quotients}
The set $H_n$ is convex-upward-closed in
$\{0,1,2\}^n\setminus\{0\}$.  For $S\in C_{2,n}$ let
$\tau_n(S)$ mean $S\subseteq H_n$.  The relation
\[
 S\equiv_n T\quad\Longleftrightarrow\quad
 \begin{cases}
  \tau_n(S)=\tau_n(T),\\
  \tau_n(SR)=\tau_n(TR)&\text{for every }R\in C_{2,n}
 \end{cases}
\]
is a congruence, and its quotient
$\mathcal R_n=C_{2,n}/{\equiv_n}$ is a finite member of
$\mathcal U_2$.  An absorption fails in $\mathcal R_n$ exactly
when a principal assignment in some multiplicative context,
including the context with no exterior factor, {has every
raw vector from the right-hand polynomial in $H_n$ and the raw vector
of the target monomial outside $H_n$}.
\end{lemma}

\begin{proof}
The set is upward closed.  Each positive point different from
$\delta_n$ has weight at least $n+1$, and
$p'_{n,i}$ has weight $2(n-1)\geq n+1$.  Hence a convex
combination followed by an upward displacement cannot produce
$\delta_n$.  If an integer point above such a combination has a
zero at $i$, every point with positive coefficient also has zero
there and must be $p'_{n,i}$.  The resulting cube point is either
$p'_{n,i}$ or a positive member of $H_n$.  A result with no zero
belongs to $H_n$ by the weight bound.  This proves closedness.

The displayed relation is a congruence by
Lemma~\ref{descent:recognition}, and its finite quotient belongs to
$\mathcal U_2$ by Proposition~\ref{cap:membership}.

Write each element of $C_{2,n}$ as a finite sum of nonzero principal
upper sets. The Frobenius identity, $(x+y)^2=x^2+y^2$, allows
squares to expand diagonally. Expanding the assigned elements and
an exterior parameter therefore expresses a capped monomial with
exponent $p$ in terms of the points $c_2(w+\sum_i p_i a_i)$.
Here each $a_i$ is nonzero, and $w=0$ denotes absence of the
exterior factor.  Closure is contained in $H_n$ exactly when all
raw points are in $H_n$.  If an absorption fails, every right raw point passes a test but
some target point fails it. Choose the principal points producing
that target. The same choices still pass in every right monomial.
Conversely, such choices give a failing assignment and test.
Variables absent from a monomial cause no difficulty, since their
choices only repeat equal summands. This proves the last assertion.
\end{proof}

For $m\geq3$ define the identity
\begin{equation}\label{clique:chi}
 \chi_m:\quad
 x_1\cdots x_m+\sum_{i=1}^m\prod_{j\ne i}x_j^2
       =\sum_{i=1}^m\prod_{j\ne i}x_j^2.
\end{equation}
Each displayed product has at least two factors.

\begin{proposition}\label{clique:independence}
For every $n,m\geq3$,
\[
 \mathcal R_n\models\chi_m\quad\Longleftrightarrow\quad n\ne m.
\]
\end{proposition}

\begin{proof}
If $n=m$, assign the principal vectors $e_i$ to the variables.
The target is $\delta_n\notin H_n$, whereas the $i$th right
point is $p'_{n,i}\in H_n$.  This refutes the identity.

Suppose $n\ne m$ and, towards a contradiction, take a principal
failure with nonzero assigned vectors $a_1,\ldots,a_m$ and
exterior vector $w$.
If the target $a=c_2(w+\sum_i a_i)$ has two
zeros, every right point retains them and fails recognition.
If it has one zero $h$ and a distinct coordinate $l$ of value one,
that one comes from the context or a unique $a_j$.  In the first
case every right point has value one at $l$ together with its zero
at $h$.  In the second case the right monomial omitting $a_j$
has a further zero at $l$.  A right point therefore fails in
either case.

The remaining bad target is $\delta_n$.  The assigned nonzero
variables then have disjoint nonempty blocks, and the context has
disjoint support.  Each right monomial omits one variable block,
so membership of all right points in $H_n$ forces every such
block to be a singleton.  Any context coordinate would have value
one in every right point, alongside its omitted zero, excluding
membership in $H_n$.  Hence there is no context support and the
$m$ singleton variable blocks exhaust all $n$ coordinates.  Thus $m=n$, a contradiction. The principal criterion now proves
the identity for every assignment and context.
\end{proof}

\begin{lemma}\label{clique:comparison}
For every $n\geq3$,
\[
 \operatorname{Id}(S_n)\subseteq\operatorname{Id}(\mathcal R_n).
\]
\end{lemma}

\begin{proof}
Take a failed absorption in $\mathcal R_n$ and choose a witness
by the principal criterion. Every right output belongs to
$H_n\subseteq K_n$.
A target with two zeros, or a positive all-one target, is already
outside $K_n$.  Otherwise the target has one zero $h$ and another
coordinate $l$ of value one.  Permute the coordinates so that $l$
immediately precedes $h$ in the cyclic order.  The set $H_n$ is
invariant under every coordinate permutation, so every right output
remains in $H_n$ and in $K_n$.  The target now has a zero with
neighbor one and is outside $K_n$.  Apply the same permutation to all assigned vectors and the exterior
vector. It commutes with capped calculations and leaves the assigned
vectors nonzero, so it gives a failure in $S_n$. The inclusion
follows by decomposing each identity into monomial absorptions.
\end{proof}

\subsection{Stabilization and comparison}

\begin{lemma}[Symmetric padding and compression]
\label{clique:compression}
For an absorption on $r$ variables with $s$ capped right monomials:
\begin{enumerate}
\item a principal failure with a coordinate positive in every right
output extends to every larger dimension;
\item a failure in $\mathcal R_n$ with $n>s$ compresses to a
failure in some $\mathcal R_d$, where
$3\leq d\leq r+s+3$, with such a common positive coordinate.
\end{enumerate}
\end{lemma}

\begin{proof}
For padding, duplicate the common positive coordinate in all assigned
vectors and in the exterior vector. This duplicates the corresponding
entry of each output. A positive right output retains a coordinate two.  A
one-zero right output has two at that coordinate, and duplication
leaves it with precisely one zero and twos everywhere else.  A
target with two zeros keeps them; an all-one target stays all one.
A target with one zero and another coordinate one keeps that
obstruction, unless its zero is duplicated, which instead gives
two zeros.  The assigned vectors remain nonzero, and a zero exterior vector
still represents omission of that factor. The duplicated coordinate
is again positive in every right output, so the operation can be repeated.

For compression, the union of the zero positions of the $s$ right
outputs has size at most $s$.  Choose a coordinate outside it.
Retain one nonzero coordinate of each assigned vector, the zero of
each one-zero right output or a two of each positive right output,
two target zero coordinates or its zero and a distinct one, and the
common positive coordinate.  An all-one target needs no additional
target coordinates.  At most $r+s+3$ coordinates have been selected.
Retain further original coordinates if necessary to reach dimension
three; the bound is at least five.

Restriction commutes with the capped calculations and leaves the
assigned vectors nonzero. Each one-zero right output keeps its zero
and has value two at every other retained coordinate, while each
positive right output keeps a coordinate two. The target remains
outside the recognized set, and the common positive coordinate remains.
Retaining further original coordinates to reach dimension three preserves
these properties. If the exterior vector becomes zero, omit the exterior
factor. This proves the compression assertion.
\end{proof}

\begin{theorem}\label{clique:stabilization}
Every identity on at most $r$ variables has the same truth value in
all $\mathcal R_n$ with
\[
 n\geq\Lambda(r):=r+3^r+2.
\]
\end{theorem}

\begin{proof}
For an absorption with $s$ right monomials put $b=r+s+3>s$.
A failure at $b$ has a common positive coordinate and pads to
every larger dimension.  A failure at any $n\geq b$ compresses
to some $d\leq b$ with a common positive coordinate, and then
pads back to $b$.  Thus its truth value is constant beyond $b$.
There are $3^r-1$ nonzero capped monomials, and equality of two
nonempty capped polynomials is equivalent to their finitely many
component absorptions in both directions.  Thus the bound $r+(3^r-1)+3=\Lambda(r)$ applies simultaneously
to all identities on at most $r$ variables.
\end{proof}

\begin{theorem}\label{clique:main}
Define
\[
 M_{\mathrm{cl}}=
 \left(\prod_{n\geq3}\mathcal R_n\right)/{\approx},
 \qquad \mathcal Z=\mathsf V(M_{\mathrm{cl}}),
\]
where $\approx$ is equality at all sufficiently large coordinates.
Then $\mathcal Z$ is a proper locally finite subvariety of
$\mathcal U_2$ with no bounded-variable identity basis.  Its
equational theory is decidable.  For an identity $\theta$ on at
most $r$ variables,
\[
 \mathcal Z\models\theta
 \quad\Longleftrightarrow\quad
 \mathcal R_{\Lambda(r)}\models\theta,
\]
and every $\mathcal R_n$, $n\geq\Lambda(r)$, satisfies all such
identities of $\mathcal Z$.
\end{theorem}

\begin{proof}
By the reduced-product criterion proved in
Theorem~\ref{eventual:main}, an identity holds in $M_{\mathrm{cl}}$
exactly when it holds in all sufficiently large $\mathcal R_n$.
Combining this with Theorem~\ref{clique:stabilization} gives the
decision rule and its simultaneous assertion. The finite quotient $\mathcal R_{\Lambda(r)}$ can be constructed
using rational linear feasibility and finitely many recognition tests.
Testing all assignments in this quotient therefore decides the identity.

Every factor belongs to $\mathcal U_2$, so its product and quotient
belong to that variety. Hence $\mathcal Z\subseteq\mathcal U_2$,
and local finiteness follows from Theorem~\ref{found:proper}.
The inclusion is proper: $\chi_3$ holds in $\mathcal Z$ but fails
in $\mathcal R_3\in\mathcal U_2$.
By Proposition~\ref{clique:independence}, every $\chi_m$ holds
in all factors except $\mathcal R_m$, so it holds in
$\mathcal Z$.  If a basis of $\mathcal Z$ used at most $r$
variables per identity, every $\mathcal R_n$ with
$n\geq\Lambda(r)$ would satisfy that basis.  It nevertheless
fails $\chi_n$, an identity of $\mathcal Z$.  This contradicts
soundness of equational deduction and proves nonfinite basability.

\end{proof}

\begin{proposition}[Comparison with the cycle reduced product]
\label{clique:strictY}
The inclusion $\mathcal Z\subsetneq\mathcal Y$ is proper.
\end{proposition}
\begin{proof}
Every identity of $\mathcal Y$ holds in all sufficiently large
$S_n$ and, by Lemma~\ref{clique:comparison}, in the corresponding
$\mathcal R_n$. It therefore holds in $M_{\mathrm{cl}}$, giving
$\mathcal Z\subseteq\mathcal Y$. Containment in $\mathcal U_2$
gives local finiteness. To prove strictness, consider
\begin{equation}\label{clique:strictidentity}
 xyzw+x^2y^2z+xyz^2w=x^2y^2z+xyz^2w.
\end{equation}
It holds in every $\mathcal R_n$.  For a bad principal target
$a=c_2(x+y+z+w+c)$, two zeros persist in both
right outputs.  With one zero and another coordinate one, the one
comes from the context or one unique variable.  The doubled-variable
sets of the right terms are $\{x,y\}$ and $\{z\}$, with empty
intersection.  Some right term consequently has value at most one
at that coordinate while retaining the zero, putting it outside
$H_n$.  For the remaining bad target $\delta_n$, all variable
blocks are disjoint and nonempty.  The first right term omits the
$w$-block.  If that block has more than one coordinate it has two
zeros; if the block is a singleton, the nonempty $z$-block supplies
a one beside its zero.  In both cases it is outside $H_n$.
The principal criterion proves \eqref{clique:strictidentity}.

In contrast it fails in every $S_n$ with $n\geq4$.  Assign
$x=e_1,w=e_2,y=e_3,z=e_4$ and use an exterior vector one at
the remaining coordinates and zero at the first four.  The target
is $\delta_n$.  The first right point has its only zero at
coordinate two and twos at its neighbors one and three, so it
belongs to $K_n$.  The second is positive with a two at coordinate
four, and also belongs to $K_n$.  Thus the principal test refutes
the identity in every such $S_n$, and the identity fails in
$\mathcal Y$.  Its validity in $\mathcal Z$ proves strictness.
\end{proof}

\begin{corollary}\label{clique:lattice}
The variety $\mathcal Z$ satisfies $G_2$ and contains the
two-element distributive-lattice variety as a proper subvariety.
\end{corollary}

\begin{proof}
We check $G_2$ directly in each $\mathcal R_n$ using the
principal criterion. A target vector of $z^2xy$ with two zeros
retains both zeros in each right output. A bad target with one zero
and a distinct coordinate one cannot obtain that one from $z$,
which is doubled. If it comes from the context, both right outputs
retain it; if it comes from $x$ or $y$, the right monomial leaving
that variable undoubled retains it. Thus at least one right output
is outside $H_n$. The remaining possible bad target, $\delta_n$,
is impossible because the nonzero vector assigned to $z$ has a
positive coordinate which its square raises to two. The opposite absorption follows from $I$. Hence $G_2$ holds
in every factor and in $\mathcal Z$.
In $C_{2,n}$ the two filters
$F_n=\mathord\uparrow(2e_1)$ and $T_n=\{(2,\ldots,2)\}$
satisfy
\[
 F_n+T_n=F_n,\quad F_nT_n=T_n,\quad
 F_n^2=F_n,\quad T_n^2=T_n.
\]
These filters form a two-element distributive lattice. Their recognition
classes are distinct: $2e_1$ has at least two zeros when $n\geq3$,
so $F_n\not\subseteq H_n$, whereas $T_n\subseteq H_n$.
The corresponding sequences remain distinct modulo eventual equality
and thus give a two-element subalgebra of $M_{\mathrm{cl}}$.
For properness, the principal filter of $\delta_n$ has false
recognition predicate, and its square $T_n$ has true predicate.
Their sequences are distinct.  Thus $M_{\mathrm{cl}}$ fails
$x=x^2$, which the lattice variety satisfies.
\end{proof}

\subsection{A criterion for squarefree absorptions}

A capped monomial is \emph{squarefree} if its nonzero exponents
are all one, and \emph{nonlinear} if an exponent is two.
For a finite nonempty variable set $U$, put
$e_U=\prod_{u\in U}u^2$.  This notation will only be used when
the product is nonempty.

\begin{theorem}[Squarefree absorption criterion]\label{clique:squarefree}
Let $q=\prod_{u\in U}u$ and let $P$ be a nonempty set of
nonlinear capped monomials supported inside $U$.  Put
$D(p)=\{u:p_u=2\}$.  The absorption $t_P+q=t_P$ fails in
$\mathcal Z$ exactly when at least one of the following holds:
\begin{enumerate}
\item $\bigcap_{p\in P}D(p)\ne\varnothing$;
\item every monomial in $P$ of proper support is precisely
$e_{U\setminus\{i\}}$ for some $i$, and the set $M$ of indices
omitted by these monomials is a proper subset of $U$.
\end{enumerate}
The case $M=\varnothing$ is allowed.  A monomial
$e_{U\setminus\{i\}}$ occurs here only when its support is nonempty.
\end{theorem}

\begin{proof}
Use the principal criterion in $\mathcal R_n$.  A target with two
zeros forces every right point to retain them.  A target with one
zero and a distinct coordinate one has that one supplied by the
context or a unique assigned variable.  To put every right point
in $H_n$, the context cannot supply it and that variable must be
doubled in every right monomial.  Conversely, if $u$ is doubled
throughout $P$, a witness exists in every dimension $n\geq3$.
Put zero at one coordinate in all vectors, put one from $u$ alone
at a second coordinate, and put two from $u$ at every remaining
coordinate.  Assign each other variable, if present, a nonzero point
at a third coordinate, and use no exterior factor.  The target
has one zero, a one at the second coordinate, and twos elsewhere.
Every right point has one zero and twos elsewhere.  This proves
exactly the first failure type.

Suppose there is no common doubled variable.  Any failure then has
target $\delta_n$.  The nonzero assigned variables occupy pairwise
disjoint nonempty blocks; the exterior vector occupies the remaining
block.  A right monomial omitting two variables has at least two
zeros.  A right monomial omitting one variable but leaving exponent
one at another gives a zero block and a nonempty block of ones.
Neither can belong to $H_n$.  Thus every right monomial of proper
support must be exactly $e_{U\setminus\{i\}}$ for an omitted
variable $i$.  For each such $i\in M$, its assigned block must
be a singleton.

If $M=U$, every variable block is a singleton.  A proper-support
monomial is present, so any context coordinate would leave a one
alongside its zero and exclude membership in $H_n$.  Hence there
is no context support and $n=|U|$ is necessary.
{Such a failure is confined to the single dimension
$|U|$, so it does not persist in the cofinite reduced product.}

If $M$ is proper, give every $M$-variable a singleton block and
distribute the remaining coordinates among nonempty blocks of the
variables in $U\setminus M$.  Use no exterior factor.  This is
possible for every $n\geq\max\{3,|U|\}$.  Each proper-support
right term has one zero and twos everywhere else.  Each full-support
right term is positive with a two, since it is nonlinear.  The
target is $\delta_n$, so the absorption fails.  When $M$ is
empty, there are no proper-support right terms, and any such partition
into variable blocks gives the same conclusion; unused coordinates
may instead be supplied by the exterior vector.

These cases exhaust principal failures in all sufficiently large
dimensions. The criterion for $\mathcal Z$ now follows from the
reduced-product criterion in Theorem~\ref{clique:main}.
\end{proof}

\section{Finite bases for proper subvarieties}
\label{clmin:section}

In this section we prove that every proper subvariety of $\mathcal Z$
has an absolute finite basis in the constant-free signature. We first
reduce an additional identity to a canonical form and then bound the
number of variables needed for the remaining identities. Put
\[
 \mathsf B_2=\C\cup\{I,P_2,G_2\},\qquad P_2:\ x^2=x^3.
\]
By Corollary~\ref{clique:lattice}, $\mathcal Z$ satisfies
$\mathsf B_2$ and contains the distributive-lattice variety
$\mathcal D$. We write $G_1$ for $xy^2=x^2y^2$, called $G$
in Section~\ref{guarded:section}. A subvariety of $\mathcal Z$
satisfying $G_1$ lies in $\mathcal G$ and is finitely based
by Theorem~\ref{guarded:hfb}.

\subsection{Normal forms and canonical additional identities}
For every finite nonempty variable set $S$ put
\[
 q_S=\prod_{i\in S}x_i,\qquad
 e_S=\prod_{i\in S}x_i^2.
\]
For $i\in S$ and $|S|\geq2$ put
\[
 h_{S,i}=x_i\prod_{j\in S\setminus\{i\}}x_j^2.
\]
All displayed products are nonempty terms. If a sum below has no
summands, we omit it, since the signature has no term for an empty sum.

\begin{lemma}[Guarded normal forms]\label{clmin:normal}
Every term is equivalent modulo $\mathsf B_2$ to a nonempty sum of
squarefree monomials and monomials of the forms $e_S$ and $h_{S,i}$.
This representation need not be unique. If one capped monomial divides
another coordinatewise, it absorbs the latter under addition.
Moreover, $\mathsf B_2$ entails
$(x+y)^2=x^2+y^2$.
\end{lemma}

\begin{proof}
By distributivity, every term is a nonempty polynomial, and $P_2$
reduces every exponent greater than two to two. Let $D$ be the nonempty
set of squared variables of a nonlinear monomial and $R$ its set of
residual variables, where $D\cap R=\varnothing$. If $R$ is nonempty, then
\begin{equation}\label{clmin:guardexpansion}
 e_D\prod_{i\in R}x_i
       =\sum_{i\in R}h_{D\cup R,i}.
\end{equation}
For $|R|=1$ the two sides coincide. For $|R|\geq2$, apply $G_2$
to two residual factors with the nonempty squared guard
$\prod_{i\in D}x_i$, and multiply by the remaining residual factors
if there are any. Each resulting monomial has one fewer residual
factor. Induction and additive idempotence give
\eqref{clmin:guardexpansion}. If $R$ is empty, the monomial is $e_D$.

If $p$ divides $t$ and they differ, write $t=pv$ with a nonempty
monomial $v$ and apply $I$.  If they coincide, apply additive
idempotence.  Finally, apply $P_2$ to $x+y$.  Its square expands
to $x^2+xy+y^2$.  Its cube, after capping and applications of
$I$, reduces to $x^2+y^2$.  Equality of the square and cube
therefore proves Frobenius.
\end{proof}

\begin{lemma}[Restriction to the target support]\label{clmin:restriction}
Let $\mathcal D\subseteq\mathcal W\subseteq\mathcal Z$, and
suppose $t_P+q=t_P$ holds in $\mathcal W$, where $q$ and the
members of $P$ are normal-form monomials.  Put $U=\supp(q)$ and
$P_U=\{p\in P:\supp(p)\subseteq U\}$.  Then $P_U$ is nonempty,
the absorption $t_{P_U}+q=t_{P_U}$ holds in $\mathcal W$, and
the latter entails the original absorption modulo $\mathsf B_2$.
\end{lemma}

\begin{proof}
In the two-element lattice algebra, assign the multiplicative identity
to the variables in $U$ and the other element to all outside variables.
If every member of $P$ involved an outside variable, this assignment
would refute the absorption. Hence $P_U\ne\varnothing$. These elements
may be assigned to variables even though neither is named in the signature.

Substitute the nonempty term $e_U$ for every outside variable.
After capping, each discarded summand becomes $e_U$. Every retained
summand divides $e_U$, so their sum absorbs it by
Lemma~\ref{clmin:normal}. The resulting identity is the restricted
absorption. Adding back the discarded summands gives the original one.
\end{proof}

\begin{lemma}[Finite-basis branches]\label{clmin:branches}
Let $\mathcal W\subseteq\mathcal Z$.
If $\mathcal W$ does not contain $\mathcal D$, it satisfies
$G_1:xy^2=x^2y^2$ and is finitely based.
Now suppose $\mathcal D\subseteq\mathcal W$ and consider an
absorption restricted as in Lemma~\ref{clmin:restriction}.
An all-squared target is absorbed modulo $\mathsf B_2$.
A nonlinear target $h_{U,i}$ is either absorbed modulo
$\mathsf B_2$ or forces $G_1$ in $\mathcal W$.
A squarefree target also forces $G_1$ if some variable is doubled in
\emph{every} right-hand summand. If a squarefree right summand occurs,
the absorption follows by divisibility, and this hypothesis cannot hold.
\end{lemma}

\begin{proof}
If $\mathcal W$ does not contain $\mathcal D$, an identity of
$\mathcal W$ failing in the two-element lattice gives a monomial
absorption $t_P+q=t_P$ failing there.  No summand of $P$ has
support contained in $U=\supp(q)$: any such support would make
the absorption valid in that lattice.  Substitute $y$ for every
variable in $U$ and $x$ for every outside variable.  The target
becomes $y^d$ for some $d\geq1$, and every right summand is
divisible by $x$.  Adding $x$ and using $I$ gives $x+y^d=x$.
When $d=1$, commutativity makes the variety trivial.  When
$d\geq2$, capping gives $x+y^2=x$.  Setting
$(x,y)=(a^2,b)$ and $(x,y)=(b^2,a)$ shows that all squares
have one common value $e$, which is an additive identity.
{Since all squares equal $e$, Frobenius and direct expansion give}
\[
 e=(a+b)^2=a^2+ab+b^2=e+ab=ab.
\]
Thus every product equals $e$, and $G_1$ holds. Since
$\mathcal W\subseteq\mathcal G$, Theorem~\ref{guarded:hfb} applies.

Every capped monomial supported inside $U$ divides $e_U$, so an
all-squared target is absorbed. If the target is $h_{U,i}$ and some
right monomial has exponent at most one at $i$, that monomial divides
the target.  Otherwise
every right monomial has exponent two at $i$.  Substitute $x$
for variable $i$ and $y$ for every other variable, and multiply
the resulting identity by $y^2$.  Capping makes its target $xy^2$
and every right summand $x^2y^2$.  Since $I$ also gives
$xy^2+x^2y^2=xy^2$, the absorption becomes $G_1$.
Here $|U|\geq2$ because the target is nonlinear.

If the target is squarefree and a variable $i$ is doubled in every
right summand, the same substitution followed by multiplication by
$y^2$ gives $G_1$. If $U=\{i\}$, take $y$ to be a fresh
variable; its square is still a nonempty term.
\end{proof}

For $U=\{1,\ldots,r\}$ with $r\geq2$, write
$E_i=e_{U\setminus\{i\}}$ and $H_i=h_{U,i}$.
For $M\subsetneq U$ define
\begin{equation}\label{clmin:K}
 K_{r,M}:\quad
 q_U+\sum_{i\in M}E_i+\sum_{i\in U\setminus M}H_i
       =\sum_{i\in M}E_i+\sum_{i\in U\setminus M}H_i.
\end{equation}
The $H$-sum is nonempty, and for $M=\varnothing$ we omit the
$E$-sums. Up to renaming variables, this identity depends only on
$r$ and $m=|M|$; we denote its renaming class by $K_{r,m}$.

\begin{lemma}[Canonical additional identities]\label{clmin:canonical}
Suppose $\mathcal W\subsetneq\mathcal Z$ contains
$\mathcal D$ and does not satisfy $G_1$.  It satisfies some
$K_{r,m}$ with $r\geq2$ and $0\leq m<r$.
More precisely, every restricted absorption of $\mathcal W$
with squarefree target and nonlinear normal-form right summands,
which fails in $\mathcal Z$, is equivalent modulo $\mathsf B_2$
to $K_{r,M}$.  Here $M$ is the union of the omitted-variable sets
of the right summands.
\end{lemma}

\begin{proof}
There is no common doubled variable, by
Lemma~\ref{clmin:branches}.  Theorem~\ref{clique:squarefree}
therefore says that every right summand with proper support is
$E_i$ for some $i$, and their omitted-index set $M$ is proper.
Every full-support summand is $e_U$ or $H_i$.
For $i\notin M$, empty intersection of the doubled-variable sets
forces $H_i$ to occur.  For $i\in M$, the term $E_i$ divides
and absorbs $H_i$ if the latter occurs.  Every member of $P$
divides $e_U$, so that summand can also be deleted when present.
The right polynomial is therefore exactly the right side of
\eqref{clmin:K}.  If $|U|=1$, every nonlinear right monomial
has a common doubled variable.  Hence $r=|U|\geq2$.

To prove existence, choose an identity of $\mathcal W$ that fails in
$\mathcal Z$. Normalize it and decompose it into monomial absorptions.
At least one absorption fails in $\mathcal Z$. Its restriction to the
target support also fails, since the restricted absorption entails the
original one. By Lemma~\ref{clmin:branches}, a nonlinear target either
gives a consequence of $\mathsf B_2$ or forces $G_1$, so neither case
can occur. A squarefree right summand would divide the squarefree target
and again make the absorption a consequence of $\mathsf B_2$. Thus the
remaining absorption has the form considered above.
\end{proof}

\begin{lemma}[Increasing rank]\label{clmin:lift}
Modulo $\mathsf B_2$, $K_{r,m}$ entails $K_{s,m}$ whenever
$s\geq r$.  The omitted-index set of size $m$ may be chosen
arbitrarily among the $s$ variables.
\end{lemma}

\begin{proof}
Choose $j\notin M$ and replace $x_j$ by $yz$, with fresh
distinct variables.  The target is the new squarefree product.
Each $E_i$ becomes the corresponding new co-square, and each
$H_i$ with $i\ne j$ becomes the corresponding new residual
monomial.  The remaining term becomes $yz e_{U\setminus\{j\}}$.
Since $r\geq2$, the square guard is nonempty, so $G_2$ expands this
term to the new $H_y+H_z$. This gives $K_{r+1,M}$, and iteration
and renaming give the result.
\end{proof}

\begin{lemma}[Enough omissions]\label{clmin:omissions}
Let $|U|=s\geq r$, and let $P$ be a nonempty set of nonlinear
capped monomials supported inside $U$, with no common doubled
variable.  Let $O$ be the union of their omitted-variable sets.
If $|O|\geq m$, then $\mathsf B_2\cup\{K_{r,m}\}$ entails
$t_P+q_U=t_P$.
\end{lemma}

\begin{proof}
For each $i\in U$, some member of $P$ has exponent at most one
at $i$ and hence divides $H_i$. Thus $t_P$ absorbs every $H_i$.
For $i\in O$, a member omitting $i$ divides $E_i$, which is
therefore also absorbed by $t_P$. Choose $m$ indices
$M\subseteq O$ and apply Lemma~\ref{clmin:lift} to obtain $K_{s,M}$.
Adding $t_P$ to this identity absorbs every right summand and gives
$t_P+q_U=t_P$. If $m=0$, no co-square is needed.
\end{proof}

\subsection{Profile compression and finite bases}

\begin{lemma}[Profile compression]\label{clmin:compression}
Suppose $t_P+q_U=t_P$ holds in $\mathcal Z$, with squarefree
target and nonlinear normal-form right summands supported inside
$U$.  Let $O$ be their omitted-variable union, suppose
$O\subsetneq U$, and put $o=|O|$.  The absorption is a substitution
instance, followed by $\mathsf B_2$ reductions, of an identity of
$\mathcal Z$ on at most
\begin{equation}\label{clmin:f}
 f(o)=o+2\bigl(2^{2^o}-1\bigr)
\end{equation}
variables.
\end{lemma}

\begin{proof}
Let $C_0=U\setminus O$. Every core variable occurs in every right
monomial, and Theorem~\ref{clique:squarefree} shows that no variable
is doubled in all of them. Hence every $v\in C_0$ is the unique
residual variable of some summand. Its profile
\[
 \mathcal H(v)=
 \{T\subseteq O:h_{U\setminus T,v}\text{ occurs in }P\}
\]
is a nonempty family of subsets of $O$, so there are at most
$2^{2^o}-1$ profiles. Partition $C_0$ according to these profiles.
Replace each class $B_0$ by $\min\{2,|B_0|\}$ fresh representatives
and leave the $O$ variables unchanged. The resulting variable set $U'$
has size at most $f(o)$.

Every original support is $U\setminus T$ for some $T\subseteq O$.
To form $P'$, replace each $e_{U\setminus T}$ by
$e_{U'\setminus T}$.  Replace a residual monomial whose residual
index belongs to $O$ by the same residual index on the new
support.  For each profile class, each of its representatives $b$,
and each $T$ in that profile, include $h_{U'\setminus T,b}$.
This gives a nonempty polynomial.

Every summand remains nonlinear. A residual $O$ variable occurs with
a squared core variable. If the original core has at least two variables,
a residual core representative occurs with another retained core
representative, which is squared in that summand. If the core has only
one variable, it was not merged, and its original square guard remains
among the $O$ variables. All-squared summands remain nonempty since
the core is nonempty.

The omitted sets $T$ are unchanged.  The doubled-variable
intersection in $P'$ is empty: each $O$ variable keeps an omission
or residual witness, and every representative has a residual
witness from its nonempty profile.  Since $O\subsetneq U$,
the validity of the original absorption in $\mathcal Z$ requires
a proper-support summand which is not a co-square, by
Theorem~\ref{clique:squarefree}.  Such a summand omits at least
two variables, or omits one and has a residual variable on its
support.  Both properties survive the construction.  Thus
$t_{P'}+q_{U'}=t_{P'}$ holds in $\mathcal Z$ by the same criterion.

Partition each original profile class $B_0$ into
$\min\{2,|B_0|\}$ nonempty blocks $B_b$, indexed by its
representatives, and substitute $\prod_{v\in B_b}v$ for $b$.
Fix the $O$ variables.  The target becomes $q_U$.
All-squared summands and residual $O$ summands recover their
original summands.  For a summand $h_{U'\setminus T,b}$,
the variables in $B_b$ become residual, and those in
$(U\setminus T)\setminus B_b$ are squared.  When $|B_b|=1$
this is already the required original summand.  When $|B_b|>1$,
the other representative in the same profile class supplies a
nonempty square guard outside $B_b$.  Formula
\eqref{clmin:guardexpansion} consequently gives exactly
\[
 \sum_{v\in B_b}h_{U\setminus T,v}.
\]
Since the profiles agree, every resulting summand occurs in $P$.
Conversely, every original core-residual summand appears in the expansion
of its own block. Thus the substituted polynomial reduces to $t_P$,
and all substituted products are nonempty.
\end{proof}

The preceding lemma applies to identities of $\mathcal Z$; the same
compression need not preserve identities of an arbitrary subvariety.
The minimal choice in the following proof ensures validity in the
ambient variety before we apply the lemma.

\begin{theorem}[Absolute finite bases for all proper subvarieties]
\label{clmin:finitebasis}
Every proper subvariety of $\mathcal Z$ is finitely based in the
constant-free signature $(2,2)$.
\end{theorem}

\begin{proof}
Let $\mathcal W\subsetneq\mathcal Z$. If it satisfies $G_1$, the
result follows from Theorem~\ref{guarded:hfb}. If it does not contain
$\mathcal D$, Lemma~\ref{clmin:branches} gives the same conclusion.
We may therefore assume that $\mathcal D\subseteq\mathcal W$
and that $G_1$ fails in $\mathcal W$.

By Lemma~\ref{clmin:canonical}, some $K_{r,m}$ holds in
$\mathcal W$. Among all such identities, choose one for which $m$
is least, and fix a corresponding $r$. Then $r\geq2$ and
$0\leq m<r$. Put
\[
 d=
 \begin{cases}
  \max\{3,r\},&m=0,\\
  \max\{3,r,f(m-1)\},&m\geq1.
 \end{cases}
\]
On a fixed alphabet of $d$ variables, let $\Sigma_d(\mathcal W)$ consist
of all identities of $\mathcal W$ between nonempty capped polynomials.
This set is finite, since there are $3^d-1$ nonzero capped monomials
and $2^{3^d-1}-1$ nonempty polynomial supports. Identities on fewer
variables are included by leaving some letters unused. We claim that
\begin{equation}\label{clmin:basis}
 \mathsf B_2\ \cup\ \{K_{r,m}\}\ \cup\
       \Sigma_d(\mathcal W)
\end{equation}
is an absolute finite basis of $\mathcal W$.
Every displayed axiom holds in $\mathcal W$.

Normalize an arbitrary identity of $\mathcal W$ and decompose it
into monomial absorptions.  By Lemma~\ref{clmin:restriction},
it is enough to derive a restricted absorption $t_P+q=t_P$ with
all right supports inside $U=\supp(q)$.
An all-squared target is absorbed.  A nonlinear residual target
is absorbed by a right monomial unless it forces $G_1$ by
Lemma~\ref{clmin:branches}; the latter alternative contradicts
our assumption.  A squarefree right summand divides a squarefree
target.  Thus only a squarefree target $q_U$ and nonlinear right
summands remain.  Their doubled-variable intersection is empty,
since otherwise $G_1$ would again hold.

Put $s=|U|$ and let $O$ be the right omitted-variable union.
If $s<r$, the absorption has at most $d$ variables and belongs
to $\Sigma_d(\mathcal W)$ after renaming.
If $s\geq r$ and $|O|\geq m$, Lemma~\ref{clmin:omissions}
derives it from $K_{r,m}$.

The remaining case has $s\geq r$ and $|O|<m$.
Then $m\geq1$ and $O\subsetneq U$, since
$|O|<m<r\leq s$.  The absorption must hold in $\mathcal Z$.
If it failed there, Lemma~\ref{clmin:canonical} would reduce it
to a $K_{s,|O|}$ valid in $\mathcal W$, contradicting minimality
of $m$.  Lemma~\ref{clmin:compression} now gives a
$\mathcal Z$-valid identity on at most
$f(|O|)\leq f(m-1)\leq d$ variables whose substitution derives
the absorption.  That compressed identity holds in
$\mathcal W\subseteq\mathcal Z$, so it is in
$\Sigma_d(\mathcal W)$ after renaming.

Adding back any discarded right summands and combining the absorptions
in both directions now derives every identity of $\mathcal W$ from
\eqref{clmin:basis}. This basis is absolute: the background axioms
$\mathsf B_2$ are finite, and the proof uses no infinite identity
theory of $\mathcal Z$ among its premises.
\end{proof}

\begin{corollary}\label{clmin:properhfb}
Every proper subvariety of $\mathcal Z$ is hereditarily finitely
based.
\end{corollary}

\begin{proof}
Every subvariety of a proper subvariety is itself proper in
$\mathcal Z$, so Theorem~\ref{clmin:finitebasis} applies to it.
\end{proof}

\begin{corollary}\label{clmin:countable}
The subvariety lattice of $\mathcal Z$ is at most countable.
\end{corollary}

\begin{proof}
There are countably many finite sets of identities in the fixed finite
signature. By Theorem~\ref{clmin:finitebasis}, each proper subvariety
is defined by one of these sets. Including $\mathcal Z$ itself still
gives an at most countable lattice; its finiteness does not follow.
\end{proof}

\begin{theorem}[An explicit limit subvariety]\label{clmin:limit}
For the explicitly defined finite recognition algebras
$\mathcal R_n$, $n\geq3$, let
\[
 M_{\mathrm{cl}}=
   \left(\prod_{n\geq3}\mathcal R_n\right)
             /{\text{eventual equality}},
 \qquad \mathcal Z=\mathsf V(M_{\mathrm{cl}}).
\]
Then $\mathcal Z$ is a limit subvariety of
$\mathsf V(\mathbb N_0,\min,+_{\mathrm{ar}})$: it is not
finitely based, and every proper subvariety is finitely based.
\end{theorem}

\begin{proof}
Theorem~\ref{clique:main} gives
$\mathcal Z\subseteq\mathcal U_2\subseteq
\mathsf V(\mathbb N_0,\min,+_{\mathrm{ar}})$ and proves that
$\mathcal Z$ has no bounded-variable identity basis.  In
particular it is not finitely based.
Theorem~\ref{clmin:finitebasis} gives an absolute finite basis
for every proper subvariety.  These are precisely the two
conditions in the definition of a limit variety.
\end{proof}

\section{Completion of the main theorem}\label{sec:conclusion}
\begin{proof}[Proof of Theorem~\ref{limit:main}]

The first assertion follows from Theorem~\ref{found:proper}.
The assertion about finite algebras is Theorem~\ref{cap:finite-members}.
Proposition~\ref{limit:joins} proves the two join equalities, and
Corollary~\ref{cap:counterexample} shows that
$\mathcal U_2$ is contained in no finite stage of either chain. The initial nonfinite-basis assertions follow from
Theorems~\ref{nfb:U2} and~\ref{nfb:L2} and
Corollary~\ref{nfb:upperinterval}. Theorem~\ref{descent:chain}
gives the next two proper nonfinitely based subvarieties. Theorems~\ref{eventual:main} and~\ref{separation:strict}
establish the first reduced-product construction, its decision
procedure and its proper containment in $\mathcal L_4$.
Theorem~\ref{smaller:main} gives the proper nonfinitely based
subvariety $\mathcal Y$. Theorem~\ref{clmin:limit} proves that
$\mathcal Z$ is a limit variety, and Proposition~\ref{clique:strictY}
shows that it is a proper subvariety of $\mathcal Y$.
Thus $\mathcal L_2$, $\mathcal L_3$ and $\mathcal L_4$ each have
a proper nonfinitely based subvariety. For every $m\geq2$ the proper inclusion
$\mathcal L_2\subsetneq\mathcal U_2\subseteq\mathcal U_m$ gives a
proper nonfinitely based subvariety, so $\mathcal U_m$ is not a limit
variety. For $m=1$, the conclusion follows from
Proposition~\ref{found:periodiccover}. The same proper nonfinitely
based subvariety shows that $\mathcal V$ is not a limit variety.
The explicitly constructed variety $\mathcal Z$ proves existence.
Since a limit subvariety must be proper in $\mathcal V$,
Theorem~\ref{found:proper} and Proposition~\ref{found:periodiccover}
show that it is periodic and locally finite. Theorem~\ref{cap:finite-members}
excludes generation by a finite algebra or a finite family of finite
algebras. Theorem~\ref{card:main} and
Proposition~\ref{smaller:upperinterval} give the cardinality assertions. Finally
Theorem~\ref{guarded:hfb} and Corollary~\ref{guarded:infinite-lattice}
give the hereditary example and its stated consequences.
\end{proof}

All proper subvarieties of $\mathcal V$ are locally finite,
but some are hereditarily finitely based and others are
nonfinitely based. The limit variety $\mathcal Z$ has at most
countably many subvarieties by Corollary~\ref{clmin:countable},
whereas the interval $[\mathcal W_\infty,\mathcal U_2]$ has
continuum many members, all nonfinitely based. The example
$\mathcal G$ shows that a hereditarily finitely based variety
may also have an infinite subvariety lattice.

%The construction establishes one explicit limit subvariety without
%classifying all such subvarieties. Three questions remain:
%\begin{enumerate}
%\item Is $\mathcal Z$ the unique limit subvariety of $\mathcal U_2$?
%\item Must every limit subvariety of $\mathcal V$ lie in $\mathcal U_2$?
%\item Can the proper subvarieties of $\mathcal Z$ be classified
%explicitly by their canonical additional identities?
%\end{enumerate}
%The power-stabilization theorem places every limit subvariety in
%some $\mathcal U_m$, but does not reduce $m$ to two. Similarly,
%the finite-basis theorem below $\mathcal Z$ does not determine
%all relations among the resulting subvarieties.

\end{document}